\documentclass[12pt]{article}

\usepackage{amssymb}
\usepackage{amsmath, amsfonts}
\usepackage[dvips]{epsfig}

\newcommand{\lbeq}[1]{{\label{OR:eq:#1}}}

\newcommand{\be}[1]{\begin{equation} \lbeq{#1}}
\newcommand{\ee}{\end{equation}}

\newtheorem{remark}{Remark}

\def\Rbb{{\mathbb R}}

\title{Chirplet approximation of band-limited, real signals made easy}

\author{J. M. Greenberg\footnote{Carnegie Mellon University, Department of Mathematical Sciences, Pittsburgh, PA  15213 (USA) {\tt greenber@andrew.cmu.edu}}
\and
Laurent Gosse\footnote{IAC--CNR ``Mauro Picone" (sezione di Bari), Via Amendola 122/D, 70126 Bari (Italy)
{\tt l.gosse@ba.iac.cnr.it}}}

\begin{document}

\maketitle

\begin{abstract}
{In this paper we present algorithms for approximating real band-limited signals by multiple Gaussian Chirps. These algorithms do not rely on matching pursuit ideas.  They are hierarchial and, at each stage, the number of terms in a given approximation depends only on the number of positive-valued maxima and negative-valued minima of a signed amplitude function characterizing part of the signal.  Like the algorithms used in \cite{gre2} and unlike previous methods, our chirplet approximations require neither a complete dictionary of chirps nor complicated multi-dimensional searches to obtain suitable choices of chirp parameters.}
\end{abstract}




\section{Introduction} 

A real-valued signal $f: \Rbb \to \Rbb$ is said to be band-limited if its Fourier transform, usually denoted by $\hat f$, has compact support. This constitutes a broad class of signals sometimes referred to as the {\it Paley-Wiener space} being of paramount importance in the applications as most of human and natural phenomena should exclude infinite frequencies. 
The classical Paley-Wiener theorem states that any band-limited signal with finite energy ({\it i.e.} such that $\int_\Rbb |f(t)|^2\ dt$ is finite) is indeed the restriction of an entire function of exponential type defined in the whole complex plane \cite{rudin}. Consequently, it would suffice to know a small part of such a signal in order to be able to extend it arbitrarily to any open set of the complex plane. However, this is unrealistic at the time being (even if some progress has been made recently, see for instance \cite{tom} and references therein) and in practice, one has still to focus on their analysis.

As the Fourier transform is supposed to be of compact support, it may seem a good idea to express numerically band-limited signals relying of optimized algorithms like the well-known {\it Fast Fourier Transform} (FFT). Especially, Shannon's sampling theorem ensures that all the signal's information can be recovered from the knowledge of a countable collection of samples. However, we shall explain in \S2.2 that this idea doesn't lead to economical (so--called {\it sparse}) representations (consult especially \cite{cj2} in this direction). An alternative can come from the decomposition into {\it chirps} which read like $A(t)\exp(j\phi(t))$; the term {\it chirplets} has been coined by Steve Mann \cite{mann} in an attempt to derive an object more sophisticated than wavelets. Briefly, if Fourier transform gives a full resolution of frequencies but a null time-resolution, wavelets offer a good compromise according to the uncertainty principle but still decompose signals onto horizontal rectangles in the time-frequency plane \cite{folle}. As chirps are endowed with a time-varying phase function $t \mapsto \phi(t)$, decomposing a signal into a superposition of chirps means that time-frequency {\it curves} are now involved. The most usual chirps are the gaussian/quadratic ones, for which both $\ln(A)$ and $\phi$ are polynomials of degree 2: we shall work in this framework hereafter.

Clearly, since chirps are more complex objects, they need more parameters to be defined correctly. Indeed, one can roughly count that 6 parameters are needed for the definition of both functions, plus 2 other parameters dealing with the time and the mean frequency around which the chirp is localized. All in all, we found at the end of \S2.3 that a signal being written as the sum of $p+1$ chirps requires at most $6p+3$ parameters to evaluate. In many cases, this constitutes a great improvement with respect to classical Fourier techniques.

The main drawback in the setup of chirp decomposition techniques is the weight of the computational effort required: according to the literature, most of the algorithms rely on matching pursuit techniques \cite{mallat} where one first considers a rather complete dictionary of chirps in order to find iteratively the ones matching the signal as best as possible (see \cite{bultan,wolves,gribon,kwok,li-stoica,li-zheng,oneill} and the recent paper \cite{candes} motivated by gravitational waves detection). However, in a very recent paper, Greenberg and co-authors \cite{gre2} proposed a completely different approach where one gets rid of the dictionary and constructs the chirp approximation by means of a simple and easy-to-implement procedure. Loosely speaking, given a complex signal in polar form $\omega \mapsto a(\omega)\exp(j\varphi(\omega))$, it amounts to seek local maxima of $a$ (we call them for instance $\omega_n$) and approximate both $\ln(a)$ and $\varphi$ by polynomials of degree 2 admitting an extremum around $\omega_n$. The procedure is first applied to $\omega_1$, the first maximum point of $a$, in order to derive 2 polynomials $A_1$ and $\phi_1$, and then it applies to the remaining signal $a(\omega)\exp(j\varphi(\omega))-A_1(\omega)\exp(j\phi_1(\omega))$ until residues become low enough. In the present paper, we propose a more elaborate algorithm than the one contained inside \cite{gre2}; it will be presented in detail within \S3. First, \S3.1 will deal with a pointwise selection procedure to compute chirps' parameters. Then in \S3.2, a mean squares $L^2$ procedure will be proposed. Both approaches will be tested on an academic example in \S4 which presents a numerical validation of these algorithms. We insist on the fact that the approach is numerically efficient and computationally cheap. Finally, \S5 will be devoted to more ``real-life" experiments: trying to decompose a signal slightly corrupted by white noise and seeking a chirp inside a stock market index.

\section{Preliminaries} 

\subsection{Specificities of real, band-limited signals} 


Our interest lies in an efficient, approximate representation of real-valued, band-limited signals $t \mapsto f(t)$.  If $f(\cdot)$ is such a time-dependent signal, it may be written in a very general way as

\begin{equation} 
\forall t \in \Rbb, \qquad f(t) = \frac{1}{\pi} \int^{\Omega}_{-\Omega} e^{jwt} {\cal H}(\omega) d \omega, \label{E1.1} 
\end{equation}

\noindent where $j=\sqrt{-1}$ is the unit of the imaginary axis, $\omega$ stands from now on for the Fourier dual variable and the function ${\mathcal H}$ rewrites like:

$$
\forall \omega \in (-\Omega, \Omega), \qquad {\cal H} (\omega) = H_e(\omega)-jH_o(\omega). \label{E1.2}
$$

\noindent The corresponding even and odd components of ${\mathcal H}$, denoted by $H_e(\cdot)$ and $H_o(\cdot)$, are smooth real-valued functions satisfying 

$$ \forall \omega \in (-\Omega,\Omega), \qquad
H_e (-\omega) = H_e(\omega)\ \ {\rm and} \ \ H_o (- \omega)=-H_o(\omega). \label{E1.3} 
$$ 

\noindent By definition of band-limited, both $H_e$ and $H_o$ vanish identically outside of the interval $(-\Omega, \Omega)$; moreover, they are assumed to satisfy 

\begin{equation}
\lim_{\epsilon \rightarrow 0^+} H_e (\Omega- \epsilon) =
\lim_{\epsilon \rightarrow 0^+} H_o (\Omega- \epsilon)
=0. \label{E1.4}
\end{equation}

\noindent The symmetries of $H_e(\cdot)$ and $H_o(\cdot)$ imply that 

\begin{equation}\forall \omega \in (-\Omega,\Omega), \qquad
\overline{\cal H} (\omega) = {\cal H}(-\omega), \label{E1.5} 
\end{equation}

\noindent which is a well-known property of the Fourier transform for real signals (the overbar stands hereafter for complex conjugate). The function ${\cal H}(\cdot)$ can be recovered from the signal via the classical Fourier transform:

\begin{equation}
{\cal H}(\omega) = \int^{\infty}_{-\infty} e^{-j\omega t}f(t) dt. \label{E1.6} 
\end{equation}

\noindent This last identify further implies that for all $\omega$,

\begin{equation}
H_e (\omega) = 2\int^{\infty}_{0} \cos(\omega t)f_e (t) dt; \qquad 
H_o (\omega) = 2 \int^{\infty}_{0} \sin(\omega t) f_o (t) dt \label{E1.8} 
\end{equation}

\noindent hold for $f_e$ and $f_o$, the even/odd parts of the signal $f$ defined for all $t \in \Rbb$ as follows:

$$
f_e(t) \stackrel{def}{=} \frac{1}{2} \Big(f(t)+f(-t)\Big)=f_e(-t) \mbox{ and } 
f_o(t) \stackrel{def}{=} \frac{1}{2} \Big(f(t)-f(-t)\Big)
=-f_o(-t). \label{E1.9}
$$

\noindent We regard (\ref{E1.6})--(\ref{E1.8}) as a sanity check on how well we are doing in synthesizing $f(\cdot)$ from ${\cal
H}(\cdot)$; that is if we employ some algorithm to compute (\ref{E1.1}), approximately, we then use that computed $f(\cdot)$ to
approximately evaluate (\ref{E1.5}) and see how well the result agrees with the original function ${\cal H}(\cdot)$.

\subsection{Standard representations of band-limited signals}

Standard representations of band-limited signals may be obtained by discretizing the integral defined in (\ref{E1.1}).  If we introduce a discrete (and finite) set of frequencies $\omega_k$

$$
\omega_k= \frac{k\Omega}{N} \ \ \ , \ \\ k \in \{-N, -N+1, ..., N\}, \label{E1.10}
$$

\noindent and exploit (\ref{E1.4}), we find that the trapezoidal rule, applied to (\ref{E1.1}), yields the approximate band-limited function, $f_N(\cdot)$, whose value reads 

$$\forall t \in \Rbb, \qquad
f_N(t) = \frac{\Omega}{N \pi} \sum^{N-1}_{k=-N+1} e^{j\frac{k \Omega t}{N}} {\cal H} \left(\frac{k \Omega}{N}\right). \label{E1.11} 
$$

\noindent The function is real-valued and satisfies $|f(t)-f_N(t)|=0 (1/N^2)$ provided the functions $H_e(\cdot)$ and $H_o(\cdot)$ are $C^2$ on $[-\Omega, \Omega]$.  More interestingly, $f_N(\cdot)$ is periodic with period, $T_N= \frac{2N\pi}{\Omega}$, and is completely determined by its values at times $t_n = \frac{n \pi}{\Omega},\ -N \leq n \leq N-1$ (this is the classical Shannon sampling theorem). To see this we note that 

$$
f_N\left(\frac{n \pi}{\Omega}\right) \stackrel{def}{=} \frac{\Omega}{N\pi} \sum^{N-1}_{k=-N+1}e^{j\frac{kn\pi}{N}}{\cal H} \left(\frac{k \Omega}{N}\right), \label{E1.12} 
$$

\noindent which yields

$$
\sum^{N-1}_{n=-N} e^{-j\frac{pn \pi}{N}} f_N\left(\frac{n\pi}{\Omega}\right) = \frac{\Omega}{N \pi} \sum^N_{k=-N+1} {\cal H} \left(\frac{k \Omega}{N}\right) \left(\sum^{N-1}_{n=-N}e^{j\frac{(k-p)n\pi}{N}}\right). \label{E1.13} 
$$

\noindent Now, one observes the following property of the exponentials:

$$
\sum^{N-1}_{n=-N} e^{j \frac{(k-p)n \pi}{N}} = \left\{ \begin{array}{llll} 
0 & , & {\rm if}\ k \neq p\\
\\
2N & , & {\rm if}\ k=p. \end{array}  \right. \label{E1.14} 
$$

\noindent The identities (\ref{E1.12})--(\ref{E1.14}) imply that for $-N+1 \leq p \leq N-1$ 

\begin{equation}
{\cal H} \left(\frac{p\Omega}{N}\right) = \frac{\pi}{2\Omega} \sum^{N-1}_{n=-N} e^{-j\frac{pn\pi}{N}}f_N \left(\frac{n \pi}{\Omega}\right) \label{E1.15} 
\end{equation}

\noindent while for $p = \pm N$, they yield: 

\begin{equation}
0 = \sum^{N-1}_{n=-N} e^{-j n \pi} f_N \left(\frac{n \pi}{\Omega}\right) = \sum^{N-1}_{n=-N} e^{jn \pi} f_N \left(\frac{n \pi}{\Omega}\right). \label{E1.16} 
\end{equation}

This identity (\ref{E1.16}) shows that if we extend (\ref{E1.15}) to the indices $p= \pm N$, then the extension is consistent with the constraint (\ref{E1.4}).  This last set of identities give an alternative means of computing ${\cal H}(\cdot)$ at the lattice points $\frac{p \Omega}{N},\ -N \leq p \leq N$.  For completeness, we record relevant identities for the coefficients $\left(H_e\left(\frac{p \Omega}{N}\right),\ H_o \left(\frac{p \Omega}{N}\right) \right)$, $0 \leq p \leq N-1$.  They are:

\begin{equation}
H_e \left(\frac{p \Omega}{N}\right) = \frac{\pi}{2 \Omega} f_N (0) + \frac{\pi}{\Omega} \sum^N_{n=1} \cos \left(\frac{pn \pi}{N}\right) f_{N,e} \left(\frac{n \pi}{\Omega}\right),\ 0 \leq p \leq N-1, \label{E1.17} 
\end{equation}

\begin{equation}
H_o \left(\frac{p\Omega}{N}\right) = \frac{\pi}{\Omega} \sum^N_{n=1} \sin \left(\frac{pn \pi}{N}\right) f_{N,o} \left(\frac{n \pi}{\Omega}\right),\ 1 \leq p \leq N-1,  \label{E1.18} 
\end{equation}

\noindent where 

$$
f_{N,e} \left(\frac{n \pi}{\Omega}\right) \stackrel{def}{=} \frac{1}{2} \left(f_N \left(\frac{n \pi}{\Omega}\right) + f_N\left(\frac{- n \pi}{\Omega}\right) \right) = f_{N,e}\left( \frac{-n \pi}{\Omega}\right),\ 1 \leq n \leq N, \label{E1.19} 
$$

\noindent and 

$$
f_{N,o} \left(\frac{n \pi}{\Omega}\right) \stackrel{def}{=} \frac{1}{2} \left(f_N
\left(\frac{n \pi}{\Omega} \right) - f_N \left(\frac{-n
\pi}{\Omega}\right) \right) = - f_{N,o} \left(\frac{-n
\pi}{\Omega}\right), \ 1 \leq n \leq N. \label{E1.20}
$$

The $\frac{2N \pi}{\Omega}$ periodicity of $f_N(\cdot)$ guarantees that $f_{N,e}\left(\frac{N\pi}{\Omega}\right) = f_N \left(\frac{-N\pi}{\Omega}\right)$ and that $f_{N,o}\left(\frac{N\pi}{\Omega}\right) =0$.  Moreover, if we extend (\ref{E1.17}) to $p=N$, then (\ref{E1.16}) guarantees that the extension is consistent with (\ref{E1.4})$_1$.  Similarly, if we extend (\ref{E1.18}) to $p=N$, we find the extension is also consistent with (\ref{E1.4}).  The standard approach outlined above will, in the limit as $N \rightarrow \infty$, yield the desired signal $f(\cdot)$ but is computationally intensive and requires a significant amount of data.  Our goal here is a more economical approximate representation of $f(\cdot)$ which, in many circumstances, requires substantially less data. 

\subsection{Chirplet representation of band-limited  signals}

We first note that if $H_e$ and $H_o$ are defined as in \S2.1 and satisfy (\ref{E1.4}), then ${\cal H}(\cdot)$ can be written in polar form (at least at points where it doesn't vanish): 

$$\forall \omega \in (-\Omega, \Omega),\qquad
{\cal H}(\omega) = A (\omega)e^{- j \psi(\omega)} \label{E1.21} 
$$ 

\noindent where the ``amplitude" is real, nonnegative and satisfies:

$$
A(\omega) \stackrel{def}{=} \left(H^2_e (\omega) + H^2_o(\omega)\right)^{1/2}=A
(- \omega), \label {E1.22} 
$$

\noindent Furthermore, $A(\omega)$ is identically zero on $|\omega|> \Omega$ and satisfies ${\displaystyle \lim_{\epsilon \rightarrow 0^+}A(\Omega- \epsilon)=0}$.  The ``phase" $\omega \mapsto \psi (\omega)$ is a smooth, odd function satisfying 

$$
\cos \psi (\omega) = \frac{H_e(\omega)}{A(\omega)}\ \ {\rm and}\ \ \sin \psi (\omega) = \frac{H_o(\omega)}{A(\omega)}. \label{E1.23} 
$$

As our {\bf canonical model} for $A(\cdot)$ we assume that it has exactly $p$ local maxima at the following distinct points, 

\begin{equation}
0 < \Omega_1 < \Omega_2 < \ldots < \Omega_p < \Omega, \label{E1.24} 
\end{equation}

\noindent and that each of these maxima is non-degenerate; i.e. $A^{(2)}(\Omega_k) < 0,\ 1 \leq k \leq p$.  The fact that
$A(-\omega)=A(\omega)$ guarantees that $A(\cdot)$ also has non-degenerate local maxima located at the symmetric set of points:

\begin{equation}
-\Omega < -\Omega_p < \ldots < - \Omega_2 < -\Omega_1 < 0. \label{E1.25} 
\end{equation}

\noindent The origin, $\omega =0$, will be either a local maximum or a local minimum of $A(\cdot)$.  In the former case we shall assume that $A^{(2)}(0)< 0$. Given the structural properties of the amplitude $A(\cdot)$ we attempt to approximate it by means of functions $A_p(\cdot)$ of the following gaussian form: 

\begin{equation}
A_p(\omega)= \alpha_0 e^{-\frac{\omega^2}{2\sigma_0}}+ \sum^p_{k=1} \alpha_k \left(e^{-\frac{(\omega-\omega_k)^2}{2 \sigma_k}} + e^{- \frac{(\omega+\omega_k)^2}{2 \sigma_k}}\right) \label{E1.26} 
\end{equation}

\noindent where $\alpha_k > 0, \sigma_k > 0$, and 

$$
0 < \omega_1 < \omega_2 < \ldots < \omega_p. \label{E1.27} 
$$

\noindent When $\omega=0$ corresponds to a minimum of $A(\cdot)$ we choose $\alpha_0=0$.  Approximate $A(\cdot)$'s of the form (\ref{E1.26}) are not band limited but they have tails which decay rapidly as $|\omega| \rightarrow \infty$ and this is adequate for our purposes. 

In \S\ref{param esti} we present two algorithms for choosing the parameters $\alpha_k > 0, \sigma_k > 0$, and the numbers $\omega_k$'s and show how these selection procedures perform on various examples. For the time being, we assume the relevant parameters are known and instead of working with ${\cal H}(\omega) = A(\omega) e^{-j\psi(\omega)}$, we replace it with 

\begin{equation}
{\cal H}_{p,1}(\omega) \stackrel{def}{=} A_p(\omega)e^{-j\psi (\omega)}. \label{E1.28} 
\end{equation}

\noindent Using this approximation to ${\cal H}(\cdot)$ we arrive at the following approximation for the signal defined by (\ref{E1.1}): 

\begin{equation}
\begin{array}{rcl}
 f_{p,1}(t) & = & {\displaystyle \alpha_0 \int^{\infty}_{-\infty} e^{- \frac{\omega^2}{2 \sigma_0}+j(\omega t-\psi(\omega))}d\omega} \\ 
\\
 &  & +{\displaystyle  \sum^p_{k=1}} \alpha_k \left(e^{j\omega_kt} 
 {\displaystyle \int^{\infty}_{- \infty}} e^{-\frac{u^2}{2 \sigma_k}+j(ut-\psi(\omega_k+u))}  du  \right.\\ 
 \\ 
 && + \left. {\displaystyle e^{-j\omega_kt} {\displaystyle \int^{\infty}_{-\infty}}e^{-\frac{u^2}{2 \sigma_k}+j(ut-\psi(-w_k+u))}du} \right)  \end{array} \label{E1.29} 
\end{equation} 

\noindent To carry this process further, we replace the terms $\psi (\pm w_k + u)$ in (\ref{E1.29}) by their local quadratic Taylor approximations around $\pm w_k$.  The result is 
$$
\psi (\omega_k + u) = \gamma_k + t_k u + \frac{\kappa_k u^2}{2},\qquad
\psi (-\omega_k + u) = - \gamma_k + t_ku- \frac{\kappa_ku^2}{2},
$$
\noindent hence
\begin{equation}
\label{E1.30}\left\{
\begin{array}{l}
\displaystyle 
\displaystyle \gamma_k = \psi (\omega_k)= - \psi (-\omega_k), \\
\displaystyle t_k = \psi'(\omega_k) = \psi'(-\omega_k), \\
\displaystyle \kappa_k = \psi'' (\omega_k)=- \psi''(-\omega_k). 
\end{array}\right.
\end{equation}

\noindent We note that this last step is equivalent to replacing ${\cal H}_{p,1}(\omega)$ with 

\begin{equation}
\begin{array}{lllll} 
{\cal H}_{p,2}(\omega) & =& \alpha_0e^{-\frac{\omega^2}{2\sigma_0}-jt_0\omega}\\
\\
&&   +{\displaystyle \sum^p_{k=1}}\alpha_k e^{-\frac{(\omega-\omega_k)^2}{2\sigma_k}-j\left(\gamma_k+t_k (\omega-\omega_k)+\frac{\kappa_k}{2}(\omega-\omega_k)^2\right)} \\
\\
&& +{\displaystyle \sum^p_{k=1}} \alpha_k 
e^{-\frac{(\omega+\omega_k)^2}{2\sigma_k}+ j \left(\gamma_k-t_k(\omega+\omega_k)+ \frac{\kappa_k}{2} (\omega+\omega_k)^2\right)}
\end{array} \label{E1.35} 
\end{equation}

\noindent and computing (\ref{E1.1})  with ${\cal H}(\cdot)$ replaced by ${\cal H}_{p,2}(\cdot)$.  The function ${\cal H}_{p,2}(\cdot)$ satisfies  $\overline{\cal H}_{p,2}(\omega)={\cal H}_{p,2}(- \omega)$ as required by (\ref{E1.5}) of ${\cal H}(\cdot)$.  The reader will note that ${\cal H}_{p,2}(\cdot)$ is merely the sum of $2p+1$ complex {\it Gaussian (or quadratic) chirps}. Exploiting the quadratic approximations defined in (\ref{E1.30}) when evaluating (\ref{E1.29}) yields, after a tedious calculation, the approximation $f_{p,2}(\cdot)$ to $f(\cdot)$:

$$
\begin{array}{rcl}
f_{p,2}(t)& = &   (2 \pi \sigma_0)^{1/2}\alpha_0e^{-\frac{\sigma_0(t-t_0)^2}{2}} \\
\\
& &+ 2^{3/2} \pi^{1/2} {\displaystyle \sum^p_{k=1}}\alpha_k 
{\displaystyle  \frac{\sigma^{1/2}_ke^{\frac{-\sigma_k(t-t_k)^2}{2(1+\kappa^2_k\sigma^2_k)}}}{(1+\sigma^2_k \kappa^2_k)^{1/4}}} \cos \left(\frac{\kappa_k\sigma^2_k(t-t_k)^2}{2(1+\kappa^2_k\sigma^2_k)} + \omega_k t-\gamma_k - \phi_k \right) \label{E1.36}
\end{array}
$$

\noindent where 

$$ 
(1+j \kappa_k\sigma_k) = (1+\kappa^2_k \sigma^2_k)^{1/2} e^{2j\phi_k} \label {E1.37} 
$$

\noindent or equivalently 

\begin{equation} \cos 2 \phi_k = \frac{1}{(1+\kappa^2_k \sigma^2_k)^{1/2}} \ \ {\rm and} \ \ \sin 2\phi_k = \frac{\kappa_k\sigma_k}{(1+\kappa^2_k\sigma^2_k)^{1/2}}. \label{E1.38} 
\end{equation}

The aforementioned approximation $f_{p,2}(\cdot)$ is simply the sum of $p+1$ real Gaussian chirps which are characterized by  $6p+3$ parameters 
\begin{itemize}
\item $(\alpha_0,\sigma_0,t_0)$, 
\item $(\alpha_k, \omega_k,\sigma_k,\gamma_k,t_k,\kappa_k)$, for $k=\{1, ...,p\}$.
\end{itemize}
Finally, $\phi_k$ is computed using (\ref{E1.38}).

\subsection{The objective of the paper} 

Approximations of band-limited signals by real-valued Gaussian chirps of the type described in (\ref{E1.36}) offers efficient representation of radar, seismic, and some biological signals.  A sampling of the literature on this subject may be found in \cite{bultan,candes,gribon,cj2,kwok,li-stoica,li-zheng,oneill,soviet} and the references contained therein.  The central problem in establishing such a formula is an efficient and natural method for obtaining the parameters characterizing the chirps.  The ``matching-pursuit" and ``ridge-pursuit" algorithms have been used by a number of authors; for details see {\it e.g.} \cite{bultan,wolves,gribon,mallat,yin}.  These typically require complicated multi-dimensional searches to obtain the chirp parameters and they also demand complete a-priori dictionary of chirps. The procedures we advance in \S\ref{param esti} to capture the parameters $(\alpha_k,\sigma_k,\omega_k)^p_{k=1}$ require no such complete dictionary and rely only on information about the signal amplitude $A(\cdot)$.  The remaining chirp parameters $(\gamma_k,t_k, \kappa_k)^p_{k=1}$ are obtained from local information about the phase, $\psi (\cdot)$, near the points $\omega_k$.  The algorithms are also hierarchical and give our approximations a multi-resolution character.  Less refined versions of these procedures were advanced in
\cite{gre2}.

\section{Parameter Selection for $A_p(\cdot)$} \label{param esti}
 
In this section we present two different algorithms for the selection of the parameters defining $A_p(\cdot)$, recall (\ref{E1.26}).  For definiteness, we assume $\omega=0$ is a local minimum of $A(\cdot)$ and thus choose $\alpha_0=0$.

\subsection{Pointwise Selection Procedure} 


We attempt to choose the parameters so that at the positive local maxima, $\Omega_1 < \Omega_2 < \ldots < \Omega_p$, one gets:
$$
A_p(\Omega_j) = A(\Omega_j),  A^{(1)}_p(\Omega_j)  =  A^{(1)}(\Omega_j)= 0,\ {\rm and} \ A^{(2)}_p (\Omega_j)  =  A^{(2)}(\Omega_j) < 0,          \eqno{(3.1)_j}  
$$ 
\noindent for $1 \leq j \leq p$. Hereafter, we denote by $g$ the Gaussian function centered in $\pm \omega_k$: 

\begin{equation} 
g(\omega; \omega_k,\sigma_k) = \left(e^{-\frac{(\omega-\omega_k)^2}{2\sigma_k}} + e^{-\frac{(\omega+\omega_k)^2}{2\sigma_k}}\right). \label{E2.2} 
\end{equation}

\noindent Then, its first derivative reads

$$
g^{(1)}(\omega, \omega_k, \sigma_k) = - \left(
\frac{(\omega-\omega_k)}{\sigma_k} e^{-\frac{(\omega-\omega_k)^2}{2\sigma_k}} + \frac{(\omega+\omega_k)}{\sigma_k} e^{-\frac{(\omega-\omega_k)^2}{2\sigma_k}}\right) \label{E2.3} 
$$

\noindent and its second derivative is

$$
g^{(2)} (\omega; \omega_k,\sigma_k)= 
\left( \left(-\frac{1}{\sigma_k} + \frac{(\omega-\omega_k)^2}{\sigma_k^2} \right) e^{-\frac{(\omega-\omega_k)^2}{2\sigma_k}}  + \left(-\frac{1}{\sigma_k} + \frac{(\omega+\omega_k)^2}{\sigma_k^2} \right) 
e^{-\frac{(\omega+\omega_k)^2}{2\sigma_k}} \right). \label{E2.4} 
$$

\noindent Hence, for $1 \leq j \leq p$, the set of equations (3.1)$_j$ rewrites as 
\begin{equation}
\alpha_j g(\omega;\omega_j,\sigma_j)=A(\Omega_j)- \sum^p_{\overset{k=1}{k \neq j}} \alpha_k g(\Omega_j; \omega_k,\sigma_k),
\label{laurent} 
\end{equation}
which boils down to: 
$$ 
\alpha_j e^{-\frac{(\Omega_j-\omega_j)^2}{2\sigma_j}} = A(\Omega_j) - \sum^p_{\overset{k=1}{k \neq j}} \alpha_k g(\Omega_j; \omega_k,\sigma_k)- \alpha_je^{-\frac{(\Omega_j+\omega_j)^2}{2\sigma_j}}. 
$$
We differentiate this equality once,
$$
\alpha_j \frac{(\Omega_j-\omega_j)}{\sigma_j} e^{-\frac{(\Omega_j-\omega_j)^2}{2 \sigma_j}} = \sum^p_{\overset{k=1}{k\neq j}} \alpha_k g^{(1)} (\Omega_j; \omega_k,\sigma_k) - \alpha_j \frac{(\Omega_j+\omega_j)}{\sigma_j}e^{-\frac{(\Omega_j+\omega_j)^2}{2\sigma_j}}, 
$$
and twice:
$$ 
\begin{array}{rcl} 
\alpha_j \left(-\frac{1}{\sigma_j}+ \frac{(\Omega_j-\omega_j)^2}{\sigma^2_j}\right) e^{-\frac{(\Omega_j-\omega_j)^2}{2 \sigma_j}} &=& A^{(2)} (\Omega_j) - {\displaystyle \sum^p_{\overset {k=1}{k\neq j}}} \alpha_k g^{(2)} (\Omega_j; \omega_k,\sigma_k)
\\

&&+ \alpha_j \left(\frac{1}{\sigma_j} - \frac{(\Omega_j+\omega_j)^2}{\sigma_j^2}\right) 
e^{ - \frac{(\Omega_j+ \omega_j)^2}{2\sigma_j}} 
\end{array}  
$$
\noindent We have in mind to propose an iterative algorithm to compute a solution to the above system of 3 equations in order to get the value of the 3 parameters; the index $n$ will refer to the successive approximations built out of this iterative process.  

\begin{itemize} 

\item For the first index $n=0$, we let 
\begin{equation} 
\alpha^0_j = A(\Omega_j),\ \omega^0_j = \Omega_j, \ {\rm and}\ \sigma^0_j = -\frac{A(\Omega_j)}{A^{(2)}(\Omega_j)},\mbox{ for } 1 \leq j \leq p. \label{E2.8} 
\end{equation} 

\item We now suppose the coefficients $\left(\sigma^n_j,\omega^n_j,\sigma^n_j\right)^p_{j=1}$ are known; we want to compute these parameters at step $n+1$. For the first index $j=1$, we let  

$$\left\{\begin{array}{lcl} 
{\cal A}^{n+1}_1 &=& \displaystyle A(\Omega_1)-\sum^p_{k=2}\alpha^n_kg \left(\Omega_1;\omega^n_k,\sigma^n_k \right)- \alpha^n_1 e^{-\frac{(\Omega_1+\omega^n_1)^2}{2 \sigma_1^n}}, \\

{\cal B}_1^{n+1} &=& \displaystyle \sum^p_{k=2} \alpha^n_k g^{(1)} \left(\Omega_1;\omega^n_k,\sigma^n_k\right) - 
\frac{\alpha^n_1(\Omega_1 + \omega^n_1)}{\sigma^n_1} 
e^{-\frac{(\Omega_1 + \omega^n_1)^2}{2\sigma_1^n}}, \\

{\cal C}^{n+1}_1 &=& \displaystyle A^{(2)}(\Omega_1) - \sum^p_{k=2}\alpha^n_kg^{(2)} (\Omega_1; \omega^n_k,\sigma^n_k)  + \alpha^n_1 \left(\frac{1}{\sigma^n_1} - \frac{(\Omega_1+\omega_1^n)^2}{(\sigma_1^n)^2}\right)e^{-\frac{(\Omega_1+\omega^n_1)^2}{2 \sigma_1^n}}
\end{array}\right.
$$

If ${\cal A}^{n+1}_1 > 0$ and ${\cal C}^{n+1}_1 < 0$, then we can compute an auxiliary quantity:

$$
{\cal D}^{n+1}_1 = {\cal B}^{n+1}_1/ {\cal A}_1^{n+1}, \label{E2.12} 
$$

\noindent and $\sigma^{n+1}_1,\omega^{n+1}_1$, and $\alpha^{n+1}_1$ are deduced as follows: 

\begin{equation} \label{E2.15} 
\left\{\begin{array}{lcl} \vspace{0.3cm} 
\sigma^{n+1}_1&=& \displaystyle  
\frac{{\cal A}^{n+1}_1}{{\cal A}^{n+1}_1({\cal D}^{n+1}_1)^2-{\cal C}^{n+1}_1}>0, \\

\vspace{0.3cm} \omega^{n+1}_1 &=& \displaystyle \Omega_1 - \sigma^{n+1}_1 {\cal D}^{n+1}_1, \\
\alpha^{n+1}_1 &=& \displaystyle {\cal A}^{n+1}_1 e^{\frac{\sigma^{n+1}_1({\cal D}^{n+1}_1)^2}{2}} > 0. 
\end{array}\right.
\end{equation} 

\item We now assume we have determined $\left(\alpha^{n+1}_j,\  \omega^{n+1}_j,\sigma^{n+1}_j\right)$ for indices $1 \leq j \leq j_0 \leq p-1$. We then fix the index $j_o$ and compute\footnote{If $j_o=p-1$, the last sums in (\ref{E2.16})-(\ref{E2.18}) are not present.} at step $j_o+1$: 

\begin{equation} 
\begin{array}{rcl}
\displaystyle {\cal A}^{n+1}_{j_{0+1}}&=& \displaystyle A\left(\Omega_{j_0+1}\right) - \sum^{j_0}_{k=1}\alpha^{n+1}_kg   \left(\Omega_{j_0+1};\omega^{n+1}_k,\sigma^{n+1}_k\right)\\
\\
 &&-\alpha^n_{j_0+1}e^{-\frac{(\Omega_{j_0+1}+ \omega^n_{j_0+1})^2}{2\sigma^n_{j_0+1}}} - 
 {\displaystyle \sum^p_{k=j_0+2}\alpha^n_k g \left(\Omega_{j_0+1}; 
 \omega^n_k,\sigma^n_k\right)}
 \end{array}  \label{E2.16} 
\end{equation}

\noindent then,

\begin{equation} 
\begin{array}{rcl} 
{\cal B}^{n+1}_{j_o+1} & = & {\displaystyle \sum^{j_o}_{k=1}} \alpha^{n+1}_k g^{(1)} \left(\Omega_{j_o+1};\omega^{n+1}_k, {\displaystyle \sigma^{n+1}_k}\right) \\
\\
&  & -\alpha^n_{j_0+1} {\displaystyle \frac{\left(\Omega_{j_o+1} + \omega^n_{j_o+1}\right)}{\sigma^n_{j_o+1}}} e^{-\frac{(\Omega_{j_o+1}+
\omega_{j_0+1})^2}{2 \sigma_{j_o+1}}} \\
\\
&  & +{\displaystyle \sum^p_{k=j_0+2}} \alpha^n_kg^{(1)} 
\left(\Omega_{j_o+1};\omega^n_k,\sigma^n_k\right), 
\end{array} \label{E2.17} 
\end{equation} 

\noindent and: 

\begin{equation}\begin{array}{rcl}
{\cal C}^{n+1}_{j_o+1} &=& A^{(2)}(\Omega_{j_o+1}) - {\displaystyle
\sum^{j_0}_{k=1}}
\alpha^{n+1}_kg^{(2)}(\Omega_{j_o+1};\omega^{n+1}_{k},\sigma^{n+1}_k)\\
\\ &&+\alpha^n_{j_0+1}{\displaystyle
\left(\frac{1}{\sigma^{n+1}_{j_o+1}}- \frac{(\Omega_{j_o+1} +
\omega^n_{j_0+1})^2}{(\sigma^n_{j_0+1})^2}\right)}
e^{-\frac{(\Omega_{j_o+1}+\omega^n_{j_0+1})^2}{2\sigma^n_{j_o+1}}}\\
\\ &&- {\displaystyle \sum^p_{k=j_o+2}} \alpha^n_k
g^{(2)}\left(\Omega_{j_o+1}; \omega^n_k, \sigma^n_2\right).
\end{array} \label{E2.18}
\end{equation} 

\noindent If ${\cal A}^{n+1}_{j_o+1}> 0$ and ${\cal C}^{n+1}_{j_o+1}<0$, then the auxiliary quantity reads 
$$
{\cal D}^{n+1}_{j_o+1}= {\cal B}^{n+1}_{j_o+1}/{\cal A}^{n+1}_{j_o+1}, 
$$

\noindent and compute $\sigma^{n+1}_{j_o+1},\ \omega^{n+1}_{j_o+1}$ and $\alpha^{n+1}_{j_o+1}$ as follows: 

\begin{equation} \label{E2.20} 
\left\{\begin{array}{lcl} \vspace{0.3cm} 
\displaystyle  \sigma^{n+1}_{j_o+1} &=& \displaystyle  \frac{{\cal A}^{n+1}_{j_o+1}}{{\cal A}^{n+1}_{j_o+1} ({\cal D}^{n+1}_{j_o+1})^2 - {\cal C}^{n+1}_{j_o+1}} > 0, \\ 
\vspace{0.3cm} \displaystyle  \omega^{n+1}_{j_o+1} &=& \displaystyle \Omega_{j_o+1} - \sigma^{n+1}_{j_o+1} {\cal D}^{n+1}_{j_o+1},\\
\displaystyle  \alpha^{n+1}_{j_o+1} &=& \displaystyle {\cal A}^{n+1}_{j_o+1} e^{\frac{\sigma^{n+1}_{j_o+1}({\cal D}^{n+1}_{j_o+1})^2}{2}} 
 > 0.
\end{array}\right.
\end{equation} 

\noindent We thus iterate until satisfactory convergence is obtained. 
\end{itemize}

\noindent The aforementioned algorithm realizes a compromise between simplicity and efficiency because it stems on inverting only half of the nonlinearities appearing in the equation (\ref{laurent}) and its derivatives in order to keep computations easy as it suffices to define the auxiliary quantity ${\cal D}$ to derive (\ref{E2.15}) and (\ref{E2.20}). This way of proceeding is related to standard iterative methods in numerical linear algebra, like for instance Jacobi or Gauss-Seidel.

\begin{remark}We note that with this method of parameter selection, the following estimate holds: $|A(\omega)-A_p(\omega)| = O (|\omega-\Omega_j|^3)$ in a neighborhood of each of the local maxima, $\Omega_j$, of $A(\cdot)$.  This pointwise selection procedure generalizes the procedure used in \cite{gre2} and captures the behavior of $A(\cdot)$ near all maxima simultaneously.
\end{remark}

 
\subsection{$L^2$ Parameter Selection} 


Again, we seek to approximate the amplitude $A(\cdot)$ by a function of the form $A_p(\cdot)$, recall (\ref{E1.26}), and choose the parameters $(\alpha_k, \omega_k, \sigma_k)$ so as to minimize the mean-squares error:

$$\begin{array}{rcl} 
||A (\cdot) -A_p (\cdot)||^2 &\stackrel{def}{=}& \displaystyle \int^{\infty}_{-\infty} |A(\omega)-A_p(\omega)|^2\ d \omega\\
\\
&=& {\displaystyle\int^{\Omega}_{-\Omega}} A^2(\omega)d \omega-2 {\displaystyle\int^{\Omega}_{-\Omega}}A(\omega)A_p(\omega) d \omega + {\displaystyle \int^{\infty}_{-\infty}} A^2_p(\omega) d \omega. \end{array} 
\label{E2.24} 
$$

\noindent For definiteness, we again assume that $\omega=0$ is a local minimum of $A(\cdot)$ and choose $\alpha_0=0$.  We also adopt the short-hand notation

\begin{equation}
G_i(\omega) = \left(e^{-{\frac{(\omega-\omega_i)^2}{2 \sigma_i}}} +
e^{-\frac{(\omega+ \omega_i)^2}{2\sigma_i}} \right)=g(\omega; \omega_i,\sigma_i), \qquad 1 \leq i \leq p, \label{E2.25}
\end{equation}

\noindent and note, for future reference, that for all $1 \leq i$ and $j \leq p$:

\begin{equation}
\begin{array}{rcl}
\left<G_i,G_j\right>= \left<G_j,G_i\right>&\stackrel{def}{=}&
{\displaystyle \int^{\infty}_{- \infty}} G_i(\omega)G_j(\omega)d
\omega\\ \\ &=& 2 \left({\displaystyle \frac{2\sigma_i\sigma_j
\pi}{\sigma_i + \sigma_j}}\right)^{1/2}
\left(e^{-\frac{(\omega_i-\omega_j)^2}{2(\sigma_i+\sigma_j)}}+
e^{-\frac{(\omega_i+\omega_j)^2}{2(\sigma_i+
\sigma_j)}}\right). \end{array} \label{E2.26}
\end{equation}

\noindent A quick calculation shows that $||A(\cdot)-A_p(\cdot)||^2$ is given by 

\begin{equation}||A(\cdot)-A_p(\cdot)||^2=
\int^{\Omega}_{-\Omega}A^2(\omega)d \omega -2\sum^p_{k=1}f_k \alpha_k
+ \sum^p_{i,j=1} \alpha_i \left<G_i,G_j\right> \alpha_j \geq 0,
\label{E2.27}
\end{equation}

\noindent where the value of $\left<G_i,G_j\right>$ has been given above and 

\begin{equation} 
f_k \stackrel{def}{=} \left<A,G_k\right>= \int^{\Omega}_{-\Omega}A(\omega)G_k(\omega) d \omega,
\qquad 1 \leq k \leq p. \label{E2.28} 
\end{equation} 

\noindent For a given choice of parameters $(\omega_i,\sigma_i)^p_{i=1}$ the $\mbox{\boldmath $\alpha$}$'s which minimize
(\ref{E2.27}) satisfy

\begin{equation} 
{\cal G} \mbox{\boldmath $\alpha$} = \mbox{\boldmath $f$}, \label{E2.29}
\end{equation}

\noindent where ${\cal G}$ is the symmetric, positive-definite, $p \times p$ matrix whose $(i,j)^{\rm th}$ entry reads: 

\begin{equation} 
{\cal G}_{i,j} = \left<G_i,G_j\right> = {\cal G}_{j,i}, \label{E2.30}  
\end{equation} 

\noindent $\mbox{\boldmath $\alpha$}$ is the $p\times 1$ column vector whose $i^{\rm th}$ entry is $\alpha_i$ and $\mbox{\boldmath $f$}$ is the $p\times 1$ column vector whose $i^{\rm th}$ entry is $f_i$.  The solution to (\ref{E2.29}) is given by $\mbox{\boldmath $\alpha$} = {\cal G}^{-1}\mbox{\boldmath $f$}$, and this implies: 

$$
||A(\cdot)-A_p(\cdot)||^2 = \int^{\Omega}_{-\Omega}A^2(\omega) d
  \omega - \mbox{\boldmath$f$}^T{\cal G}^{-1} \mbox{\boldmath $f$}
  \geq 0. \label{E2.32}
$$ 

\noindent Thus, the task before us is to choose $(\omega_i,\sigma_i)^p_{i=1}$ to maximize $\mbox{\boldmath $f$}^T {\cal G}^{-1}\mbox{\boldmath $f$}$ where again ${\cal G}^{-1}$ is the symmetric, positive-definite inverse of ${\cal G}$ defined in
(\ref{E2.30}) above. At first blush, the problem of maximizing $Q=\mbox{\boldmath $f$}^T{\cal G}^{-1}\mbox{\boldmath $f$}$ looks rather daunting but, in fact, this problem has a rather simple structure. 

\begin{itemize} 
\item For any integer $k$ between 1 and $p$ we let $\beta_k$ be one of the two parameters $\omega_k$ or $\sigma_k$ and we observe that:

\begin{equation}
\frac{\partial Q}{\partial \beta_k} = \mbox{\boldmath $f$}^T{\cal
G}^{-1},_{\beta_k}\mbox{\boldmath $f$} +2 \mbox{\boldmath $f$}^T{\cal
G}^{-1} \mbox{\boldmath $f$},_{\beta_k}. \label{E2.33}
\end{equation} 

\noindent Noting that $\mbox{\boldmath $f$}^T{\cal G}^{-1}=\mbox{\boldmath $\alpha$}^T$ and that $\mbox{\boldmath $f$},_{\beta_k}$ is the $p \times 1$ column vector whose $k^{\rm th}$ component is $f_{k,\beta_k}= \left<A,G_{k,\beta_k}\right>$ and other components zero, we see that the second term on the right-hand side of (\ref{E2.33}) is $2 \alpha_kf_{k,\beta_k}$.  Moreover, if we exploit the identity,

\begin{equation} 
{\cal G}^{-1},_{\beta_k} =-{\cal G}^{-1}{\cal G},_{\beta_k} {\cal G}^{-1}, \label{E2.34} 
\end{equation}

\noindent we find that the first term on the right-hand side of (\ref{E2.33}) is 
$$
-2\alpha_k {\displaystyle \sum^p_{\overset{j=1}{j\neq k}}} \left<G_k,G_j\right>\!_{\underset{,\beta_k}{\ \ }}\alpha_j
-\alpha^2_k\left<{G}_k,{G}_k\right>\!\!_{\underset {,\beta_k}{\ \ }}
$$
and again the $\mbox{\boldmath $\alpha$}$'s are the solution of (\ref{E2.29}).  These observations imply that

$$ 
\frac{\partial Q}{\partial \beta_k} = 2 \alpha_k \left(\frac{\partial
f_k}{\partial \beta_k} - \sum^p_{\overset{j=1}{j \neq k}}
\frac{\partial \left<G_k,G_j\right>}{\partial \beta_k} \alpha_j\right)
- \alpha^2_k \left<G_k,G_k\right>\!\!_{\underset{,\beta_k}{\ \ }} \label{E2.35}
$$ 

\noindent where again ${\cal G} \mbox{\boldmath $\alpha$} = \mbox{\boldmath $f$}$. Given the particularly simple structure of $\frac{\partial Q}{\partial \beta_i}$ we solve the problem of maximizing $Q$ by a ``steepest ascent" type algorithm as we explain now.

\item So we assume now $(\omega^n_i,\sigma^n_i)^p_{i=1}$ being known.  With these data, and the explicit formula (\ref{E2.25}) for $G_i(\omega)$, we explicitly compute the $\left<G_i,G_j\right>$'s and $\frac{\partial}{\partial \beta_i} \left<G_i,G_j\right>$'s at the $n^{\rm th}$ data set.  Once again $\beta_i=\omega_i$ or $\sigma_i$.  We denote the results by $\left<G_i,G_j\right>^n$ and $\frac{\partial}{\partial \beta_i}\left<G_i,G_j\right>^n$ respectively.  We sample the amplitude $A(\cdot)$ at points $$\omega_p=\frac{p\Omega}{N} \ \ \ , \ \ \ -N+1\leq p \leq N-1,$$ and do the same for the functions $G_i(\cdot)$ and $\frac{\partial G_i}{\partial \beta_i}(\cdot)$.  These are of course also evaluated with the $n^{\rm th}$ data set and the results are superscripted with the index $n$.  We then derive approximate values for $f_i$ and $\frac{\partial f_i}{\partial \beta_i}$ using the discrete inner products defined below to replace the integrals (see (\ref{E2.28}))

$$
f^n_i \approx \frac{\Omega}{N}  \sum^{N+1}_{p=-N+1} G^n_i \left(\frac{p \Omega}{N}\right)A\left(\frac{p\Omega}{N}\right) \label{E2.37} 
$$

\noindent and 

$$ 
\frac{\partial f^n_i}{\partial \beta_i} \approx \frac{\Omega}{N} \sum^{N-1}_{p=-N+1} \frac{\partial G^n_i}{\partial \beta_i} \left(\frac{p\Omega}{N}\right)A\left(\frac{p\Omega}{N}\right). \label{E2.38} 
$$

\noindent Next we obtain $\mbox{\boldmath $\alpha$}^n$ by solving ${\cal G}^n \mbox{\boldmath $\alpha$}^n= \mbox{\boldmath $f$}^n $ in order to compute,

\begin{equation} 
\frac{\partial Q^n}{\partial \beta_i} = 2 \alpha^n_i
\left(\frac{\partial f^n_i}{\partial \beta_i} -
\sum^p_{\overset{j=1}{j \neq 1}} \frac{\partial \left<G_i,G_j\right>^n
}{\partial \beta_i} \alpha^n_j\right) - (\alpha^n_i)^2
\left<G_i,G_i\right>^n\!\!_{\underset{,\beta_i}{\ \ }}, \label{E2.40}
\end{equation} 

\noindent and use (\ref{E2.40}) to update the parameters as follows: 

$$
\omega^{n+1}_i = \omega^n_i + \Delta \frac{\partial Q^n}{\partial \omega_i}, \qquad
\sigma^{n+1}_i = \sigma^n_i + \Delta \frac{\partial Q^n}{\partial \sigma_i}. 
$$ 

\noindent Unless $\frac{\partial Q^n}{\partial \omega_i} = \frac{\partial Q^n}{\partial \sigma_i} = 0, \ 1 \leq i \leq p, \ \ Q$ will increase for $\Delta$ small enough. We iterate the aforementioned process till convergence.
\end{itemize}

\section{A first ``academic" numerical example} 

\subsection{First test} 

To see how these selection processes perform we consider the following example.  We let 

\begin{equation} 
A(\omega) = \left\{\begin{array}{ll} 
0, & - \infty < \omega \leq -2\\
(4-\omega^2)^2\left({\displaystyle \frac{1}{2}} + \omega^2\right), & \ -2 \leq \omega \leq 2\\
0, & 2 \leq \omega < \infty. \end{array} \right.  \label{E2.43} 
\end{equation} 

\noindent The origin $\omega =0$ is a local minimum of $A$ and the points $\omega = \pm 1$ are the maxima of $A(\cdot)$ with $A(\pm1) = 13.5$.  The second derivative at these points is $A^{(2)}{(\pm 1)} = -36$.  Our approximations are of the form:

\begin{equation} 
A_1(\omega)= \alpha_1 \left(e^{-\frac{(\omega-\omega_1)^2}{2\sigma_1}} + 
e^{- \frac{(\omega+\omega_1)^2}{2\sigma_1}}\right), \qquad \omega \in \Rbb.         \label{E2.44}
\end{equation}

\begin{figure}[ht]
\centerline{\epsfig{file=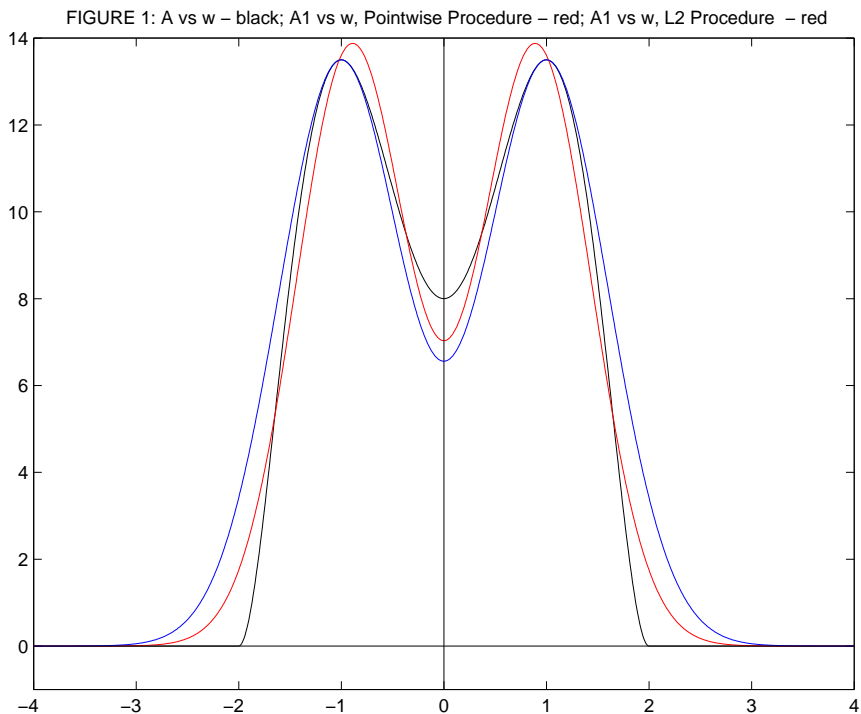,angle=0,width=4.5in,height=3.75in}}
\caption{Comparison between pointwise and $L^2$ selection procedures with example (\ref{E2.43})--(\ref{E2.44}): original signal $A$ is in black, pointwise algorithm in blue and $L^2$ in red.}
\label{jim1}
\end{figure}

\noindent The pointwise selection procedure yielded the following results: $\alpha_1 = 13.4515$, $\omega_1=1.0074$ and $\sigma_1 = 0.3595$, which in turn yields $\underset{\omega}{\rm max}A_1(\omega) = A_1(\pm1)=13.5$.  The $L^2$ procedure furnished $\alpha_1=13.8189$, $\omega_1=.8974$ and $\sigma_1 = 0.2942$; in this case $\underset{\omega}{\rm max}A_1(\omega) = A_1 (\pm .9)=13.8782$.  The maximal value of $Q$ is 390.9413 and the $L^2$ norm of $A(\cdot)$ is 395.0055. Figure \ref{jim1}, below, shows a graph of $A(\cdot)$ along with graphs of $A_1(\cdot)$ for each selection procedure.

The results for this simple example are not atypical of the general case; namely, one application of either selection process typically yields an approximation, $A_p(\cdot)$, which is qualitatively similar to the given amplitude, $A(\cdot)$, but may differ quantitatively from $A(\cdot)$. This defect may be overcome by repeated applications of the procedure.

\subsection{Hierarchical Refinements of Selection Procedures} 

We assume we have already applied either the pointwise or $L^2$ selection procedures $n-1$ times and constructed the functions
$A_{N_k,k}(\cdot) , \ 1 \leq k\leq n-1$, each which are sums of $N_k$ Gaussians.  We let

$$ \forall \omega \in (-\Omega,\Omega), \qquad
A_n(\omega)=A_0(\omega)- \sum^{n-1}_{k=1}A_{p_k,k}(\omega), \label{E2.45} 
$$

\noindent where $\omega \mapsto A_o(\omega)$ is the original amplitude, previously denoted by $A(\cdot)$.  As defined $A_n(\cdot)$ is even and rapidly decreasing as $|\omega|\rightarrow \infty$.  The principal qualitative difference between $A_n(\cdot)$ and $A_o(\cdot)$ is that $A_n(\cdot)$ takes on both positive and negative values.  For definiteness, we assume now that $\omega=0$ is a maximum of $A_n(\cdot)$ and that the typical structure of $A_n(\cdot)$ is as follows:

\begin{itemize} 
\item there are exactly $q_n$ negative-valued, local minima of $A_n(\cdot)$ at points 

$$ 
0 <\underline{\Omega}^n_1 < \underline{\Omega}^n_2 < \ldots < 
\underline{\Omega}^n_{q_n} \label{E2.46} 
$$

\noindent and each minima is non-degenerate; 

\item there are exactly $p_n$ positive-valued, local maxima of $A_n(\cdot)$ at points 

$$
0 < \overline{\Omega}^n_1 < \overline{\Omega}^n_2 < \ldots < \overline{\Omega}^n_{p_n} \label{E2.47} 
$$

\noindent and each maxima is non-degenerate; 

\item the minima and maxima are not necessarily interlaced since $A_n(\cdot)$ may have positive valued local minima and negative valued local maxima.
\end{itemize} 

\begin{figure}[ht]
\centerline{\epsfig{file=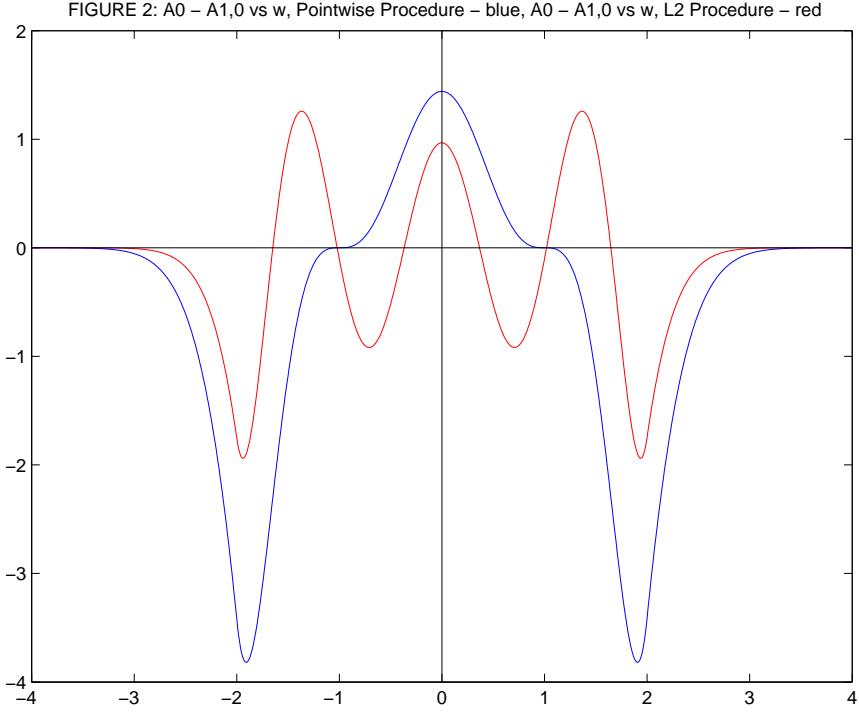,angle=0,width=4.5in,height=3.75in}}
\caption{Difference $A_0(\omega)-A_{1,0}(\omega)$: pointwise selection in blue, $L^2$ selection in red.}
\label{jim2}
\end{figure}

\noindent Graphs of $A_1(\cdot) = A_0(\cdot)-A_{1,0}(\cdot)$ for the two selection procedures used in our previous example are shown in Figure \ref{jim2}. 
Our induction step is to replace $A_n(\cdot)$ by a sum of $N_n=q_n+p_n+1$ even Gaussians:  

\begin{equation} 
\begin{array}{rcl}
A_{N_n,n}(\omega) & = & \alpha^n_0
e^{-\frac{\omega^2}{2\overline{\sigma}_0^n}} + {\displaystyle
\sum^{p_n}_{k=1}} \alpha^n_k
\left(e^{-\frac{(\omega-\overline{\omega}^n_k)^2}{2\overline{\sigma}^n_k}}
+ e^{-\frac{(\omega+\overline{\omega}^n_k)^2}{2
\overline{\sigma}^n_k}}\right) \\ 
&& - {\displaystyle
\sum^{q_n}_{k=1}}\beta^n_k
\left(e^{-\frac{(\omega-\underline{\omega}^n_k)^2}{2
\underline{\sigma}^n_k}}+ e^{-\frac{(\omega+
\underline{\omega}^n_k)^2}{2 \underline{\sigma}^n_k}}\right)
\end{array} \label{E2.48}
\end{equation} 

\noindent where the $\alpha$'s, $\beta$'s, \underline{$\sigma$}'s, and 
$\overline{\sigma}$'s are positive and the \underline{$\omega$}'s and $\overline{\omega}$'s satisfy 

$$
0 \leq \underline{\omega}^n_1 < \underline{\omega}^n_2 < \ldots < \underline{\omega}^n_{q_n}, \qquad
0 < \overline{\omega}^n_1 < \overline{\omega}^n_2 < \ldots < \overline{\omega}^n_{p_n}. \label{E2.50} 
$$

The pointwise selection procedure generates the unknown parameters by insisting that $A_{N_n,n}(\cdot),\ A^{(1)}_{N_n,n}(\cdot)$, and $A^{(2)}_{N_n,n}(\cdot)$ match $A_n(\cdot), \ A^{(1)}_n(\cdot)$, and $A^{(2)}_n (\cdot)$ at the local maxima $\{0,
\{\overline{\Omega}^n_k\}^{p_n}_{k=1}\}$ and local minima $\{\underline{\Omega}^n_k\}^{q_n}_{k=1}$. The $L^2$ selection procedure chooses the coefficients to minimize $||A_n(\cdot)-A_{N_n,n}(\cdot)||^2$.  This problem has the same structure as the optimization problem discussed in detail earlier, and may be solved by the ``steepest-ascent" algorithm.  The optimal solution satisfies

\begin{equation} 
0 \leq ||A_{n+1}(\cdot)||^2 = ||A_n(\cdot)-A_{N_n,n}(\cdot)||^2 = ||A_n(\cdot)||^2-Q_{n,{\rm max}}, \label{E2.51} 
\end{equation} 

\noindent where $(\underline{\omega}^n,\underline{\sigma}^n; \overline{\omega}^n, \overline{\sigma}^n) \rightarrow Q_n
=||A_{N_{n,n}}(\cdot)||^2$ when the $\mbox{\boldmath $\alpha$}$'s and $\mbox{\boldmath $\beta$}$'s are the least squares parameters corresponding to the given choice $\underline{\omega}^n, \underline{\sigma}^n,\overline{\omega}^n$, $\overline{\sigma}^n$ and $Q_{n,{\rm max}} = \underset{(\underline{\omega}^n,\underline{\sigma}^n,\overline{\omega}^n,\overline{\sigma}^n)}{\rm max} Q_n$. 
The equality (\ref{E2.51}) implies that 

$$ 
0 \leq ||A_{n+1}(\cdot)||^2 =
||A_0(\cdot)-\sum^n_{j=0}A_{N_j,j}(\cdot)||^2=||A_0(\cdot)||^2-\sum^n_{j=0}Q_{j,{\rm
max}}, \label{E2.52}
$$

\noindent and provides us with a stopping criteria for the number of applications of the $L^2$ selection procedure.  Specifically, we stop when $||A_{n+1}(\cdot)||^2 \leq \epsilon_{\rm stop}$, a preassigned tolerance. The stopping criteria for the pointwise selection procedure is equally easy; we stop when 

\begin{equation} 
\underset{\omega}{\rm max} (A_{n+1}(\omega)) \leq \epsilon_{\rm stop}\ \ {\rm
and} \ \ \underset{\omega}{\rm min} (A_{n+1}(\omega)) \geq -
\epsilon_{\rm stop}. \label{E2.53}
\end{equation} 

\begin{figure}[ht]
\centerline{\epsfig{file=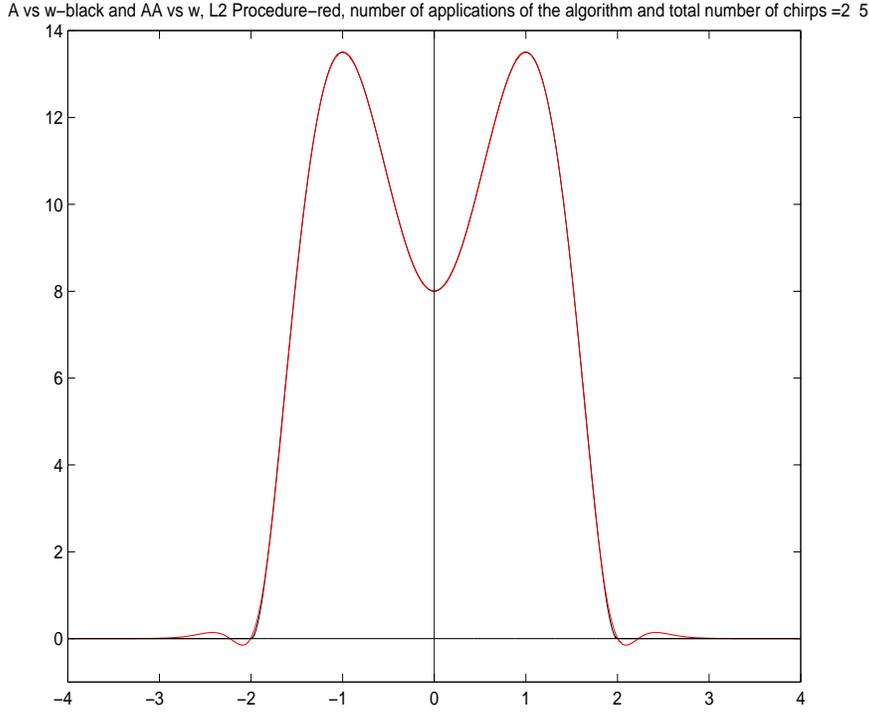,angle=0,width=4.5in,height=3.75in}}
\caption{$L^2$ hierarchical refinement procedure (to be continued in Fig. \ref{jim3'}): Original signal $A$ is in black and $L^2$ in red.}
\label{jim3}
\end{figure}

\begin{figure}[ht]
\centerline{\epsfig{file=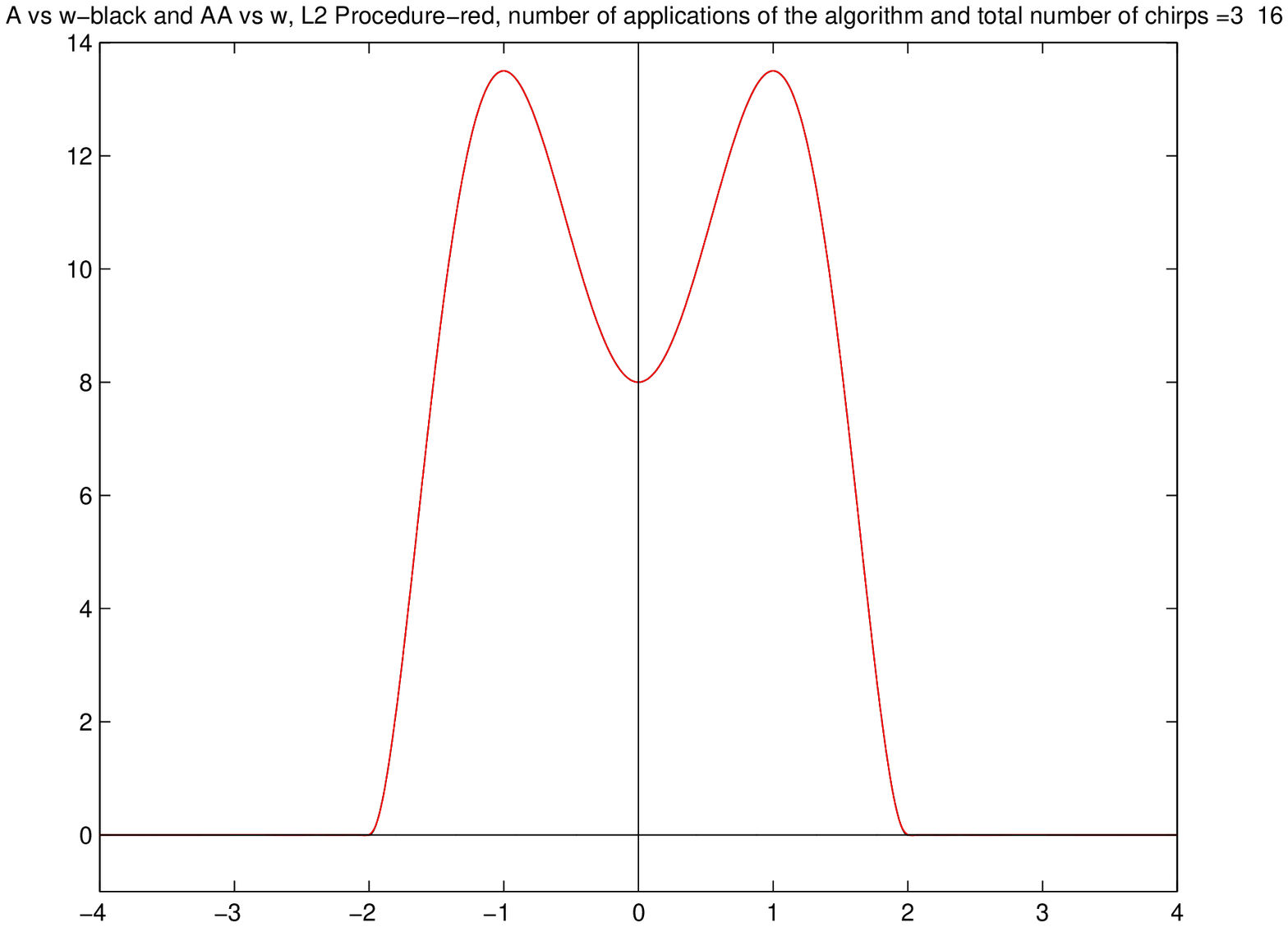,angle=0,width=4.5in,height=3.75in}}
\centerline{\epsfig{file=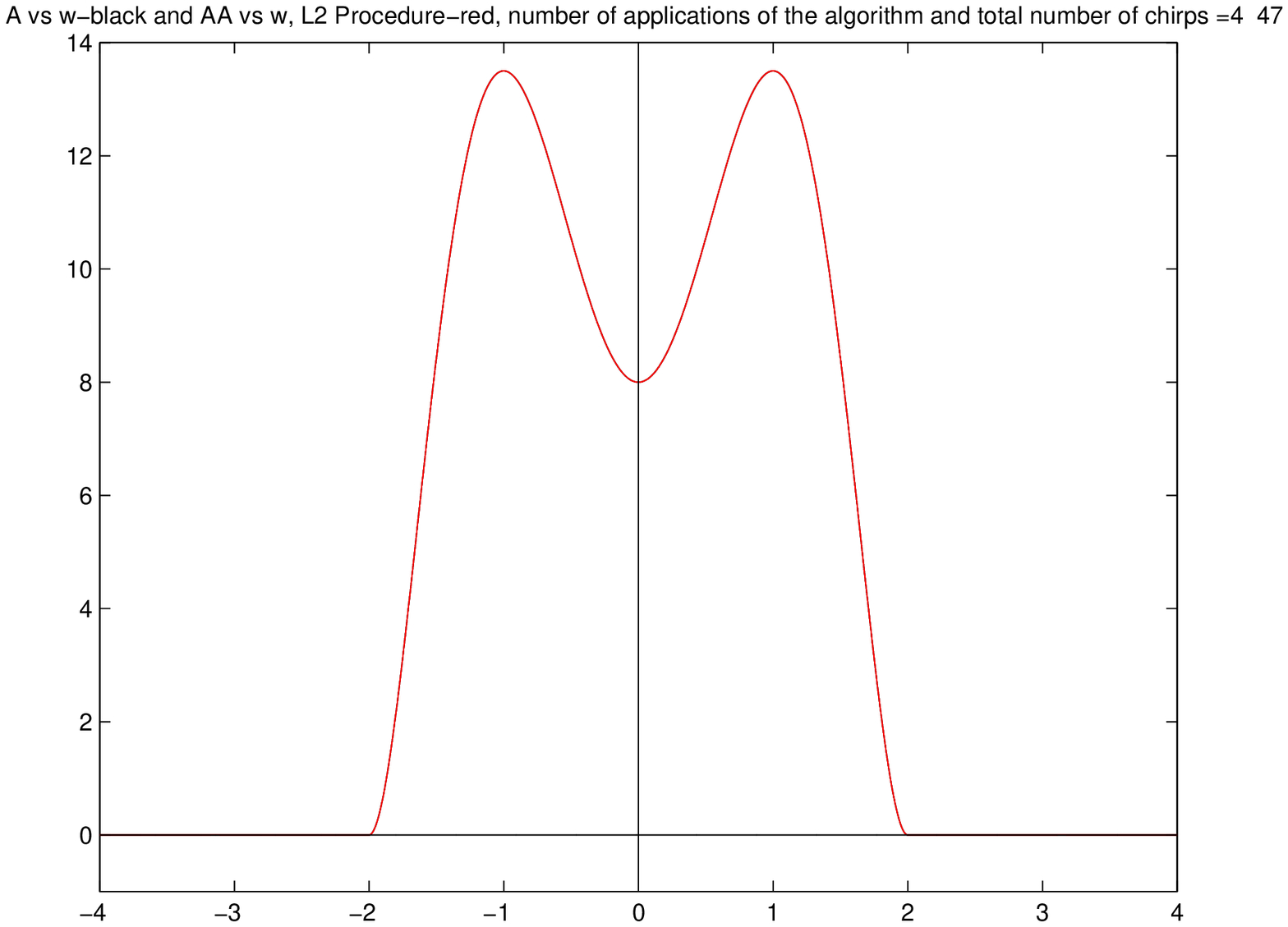,angle=0,width=4.5in,height=3.75in}}
\caption{$L^2$ hierarchical refinement procedure, continued from Fig. \ref{jim3}.}
\label{jim3'}
\end{figure}

Figures \ref{jim3} and \ref{jim3'} shows the results of applying the hierarchical $L_2$ procedure 2,3, and 4 times, respectively.  The ``black'' curves on each figure are the original amplitude, $A(\cdot)$, and the ``red''curves are the composite hierarchial gaussian 
approximations.  In applying the $L^2$ procedure twice we went from 1 to 5 even gaussians; in the third application we went from 5 to 16 even gaussians, and in the fourth application from 16 to 47 even gaussians.  After four applications of the algorithm, the curves become indistinguishable.


\section{Two ``real-life" applications} 

\subsection{Chirp decomposition of a noisy signal} 

A band-limited signal slightly corrupted by white noise isn't band limited anymore; however, if we set up the preceding algorithms with a value of $\Omega$ being large enough, it may be possible to recover a correct approximation. First of all, we set up the following data: consider the even amplitude function displaying only two bumps,
\begin{equation} 
A(\omega)=\frac{\exp(-a|\omega|^3)-\exp(-b|\omega|^3)}{b-a}; \qquad a=0.8,\ b=1.3.
\label{lolo}
\end{equation}
We couple it with different phase functions of increasing complexity: first, we chosed a cubic one, $\phi(\omega)=\frac{\omega^3}{50}$ and then $\phi(\omega)=\pi(1-\exp(-\omega^2))\sin(2\omega)$. The discretization grid is 512 points from $t=-5.1$ to $t=5.12$. We don't insist on the fact that in case the original phase function $\phi$ is a polynomial of degree 2, its recovery is exact.
\begin{figure}[ht]
\centerline{\epsfig{file=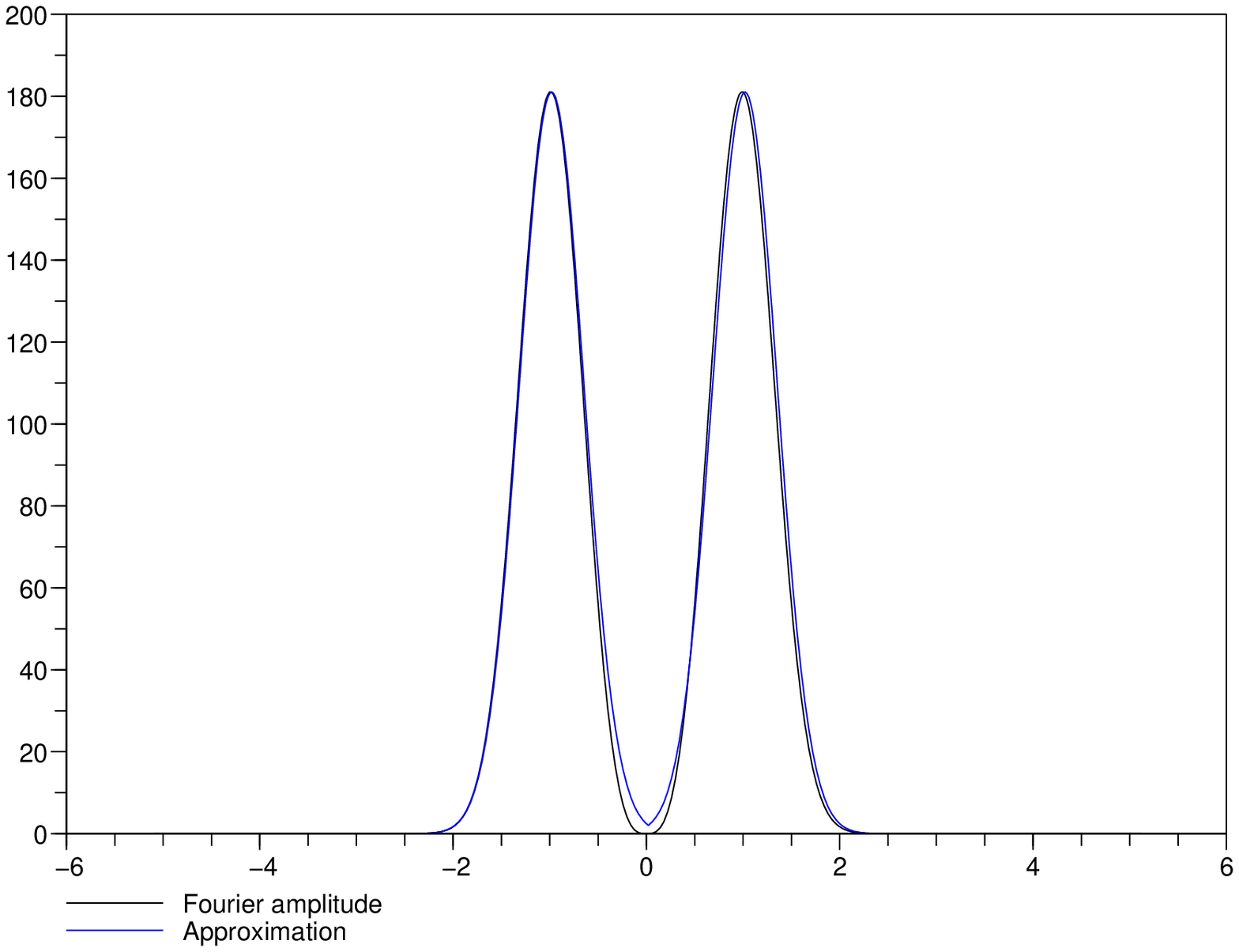,angle=0,width=2.5in,height=3.75in}
\epsfig{file=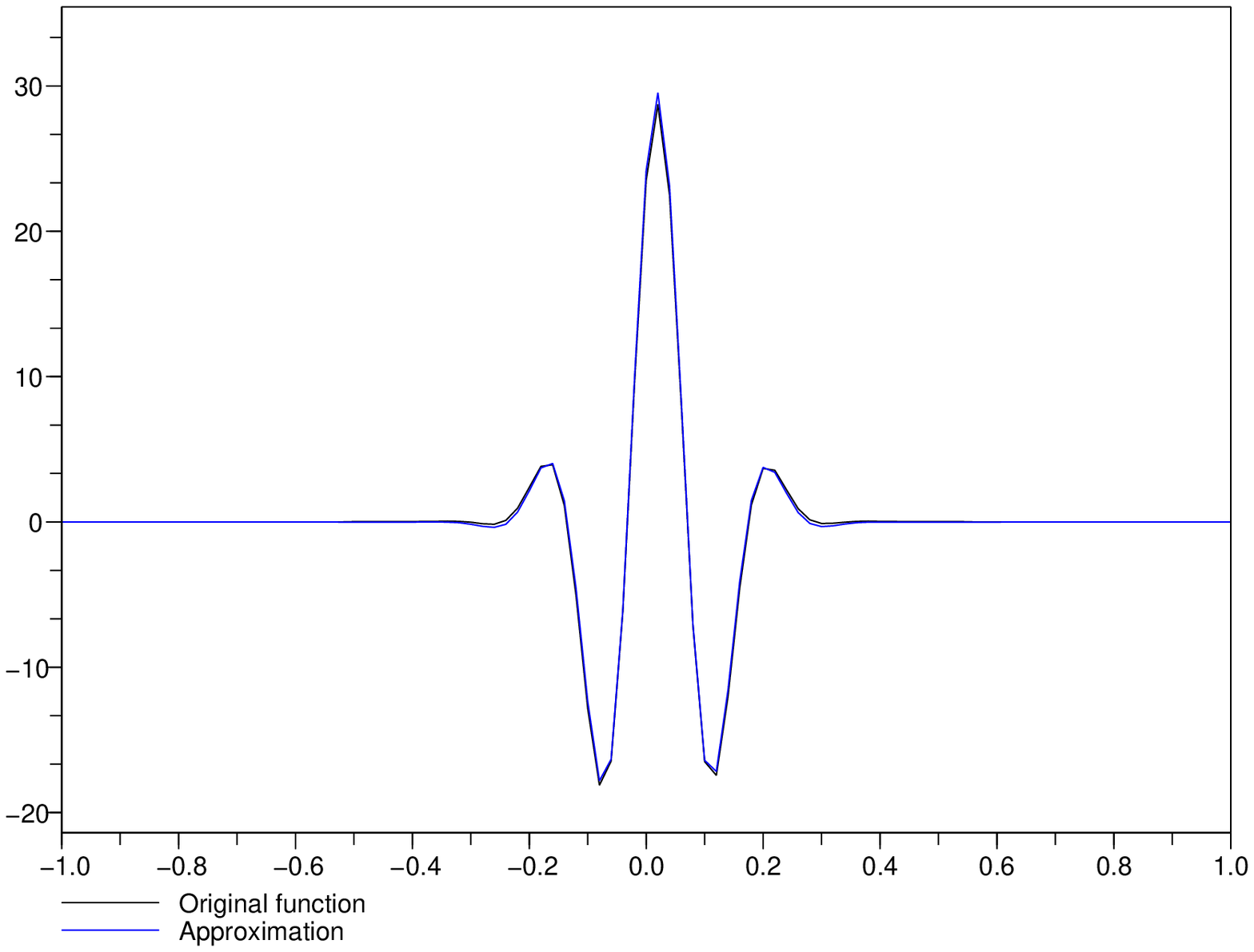,angle=0,width=2.5in,height=3.75in}}
\centerline{\epsfig{file=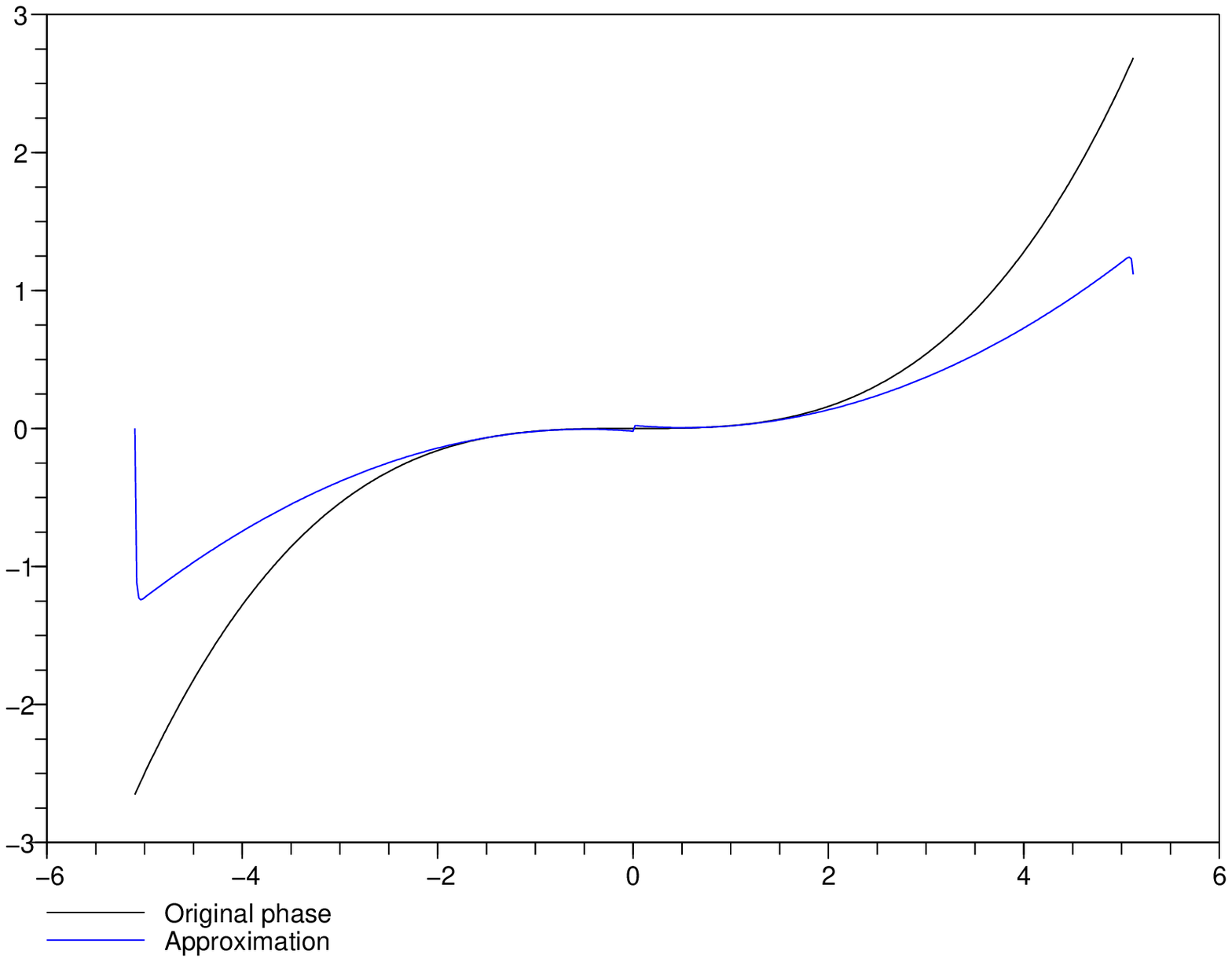,angle=0,width=2.5in,height=3.75in}
\epsfig{file=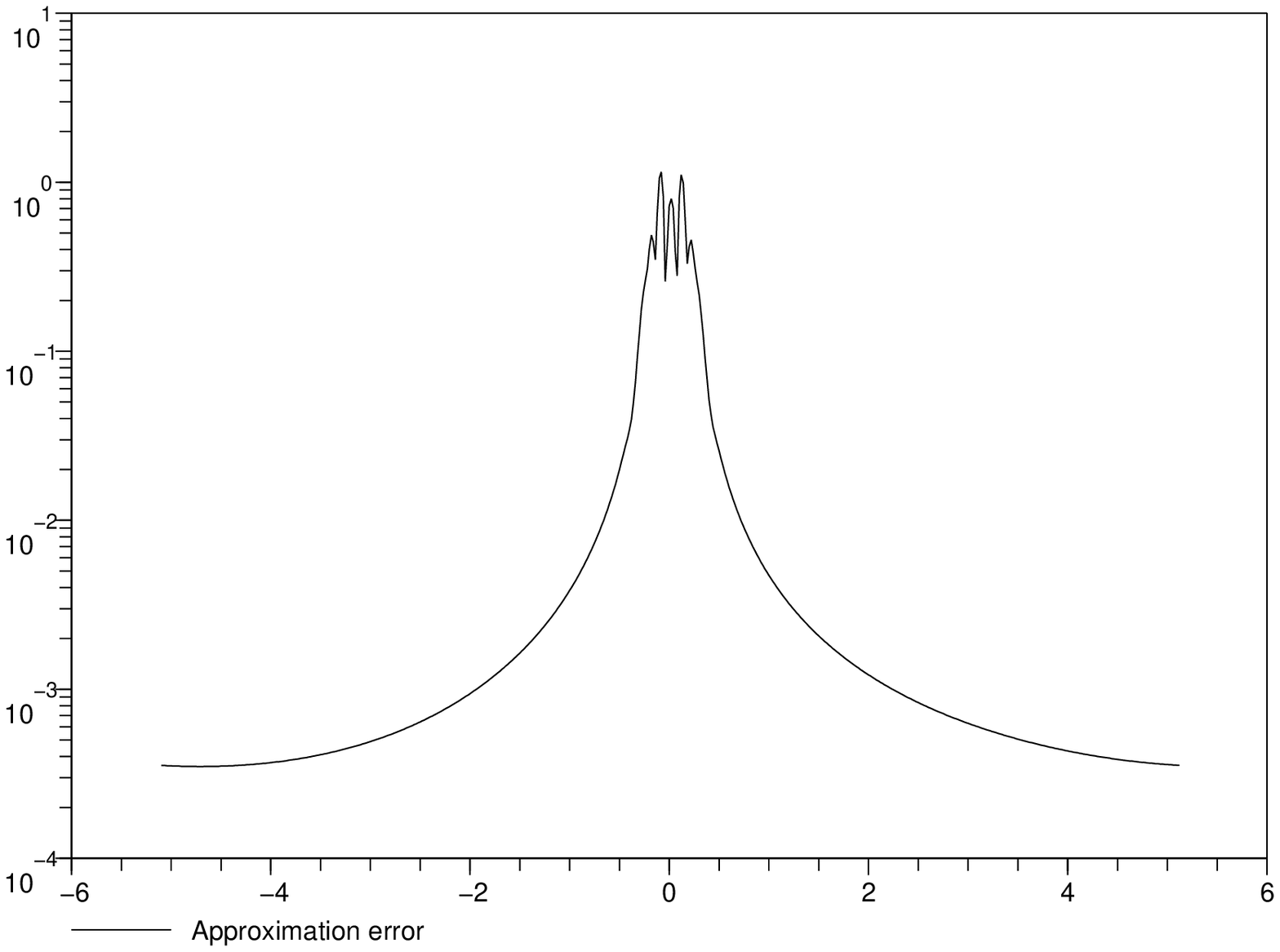,angle=0,width=2.5in,height=3.75in}}
\caption{Pointwise procedure selection on example (\ref{lolo}) for cubic phase.}
\label{cubic}
\end{figure}
\begin{figure}[ht]
\centerline{\epsfig{file=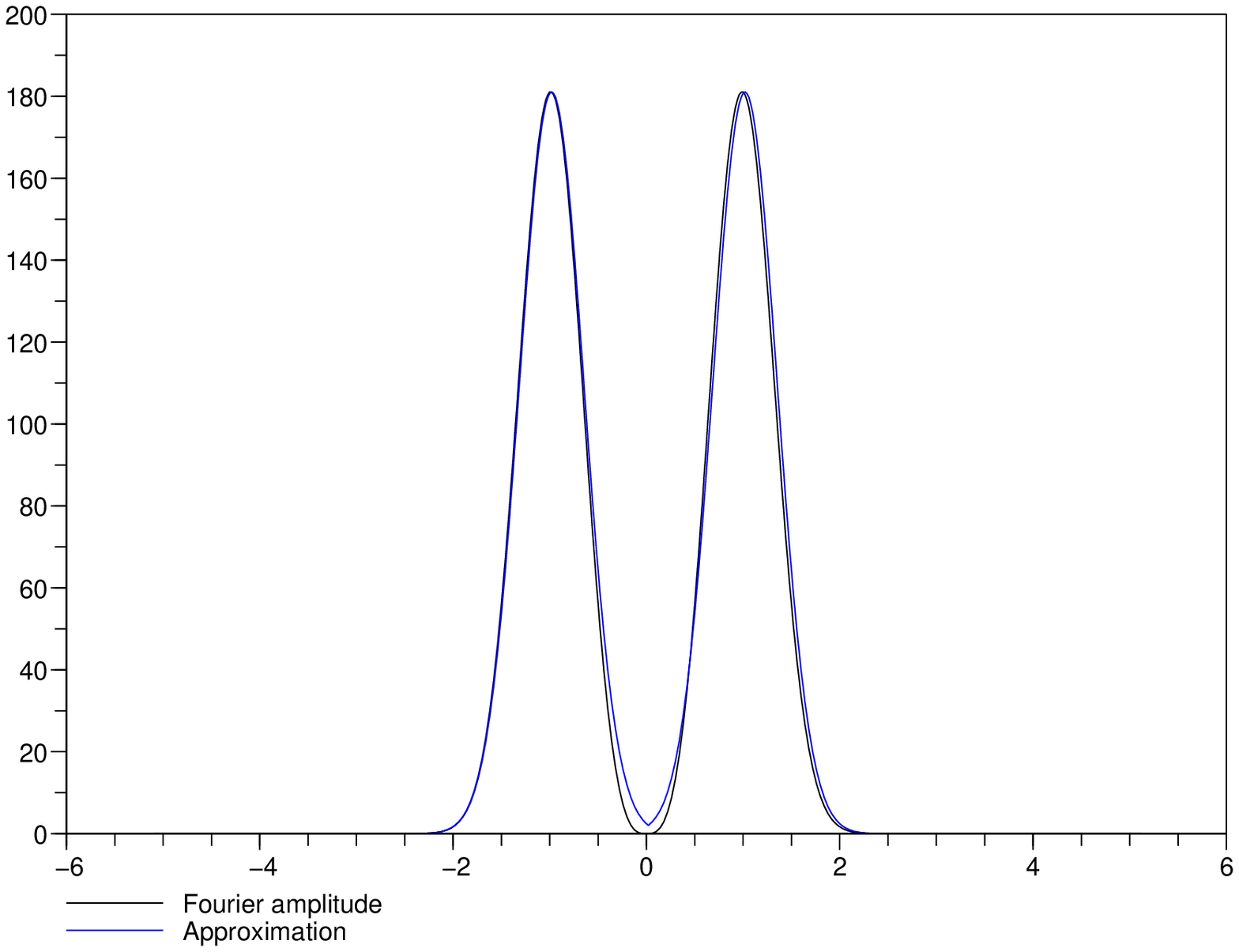,angle=0,width=2.5in,height=3.75in}
\epsfig{file=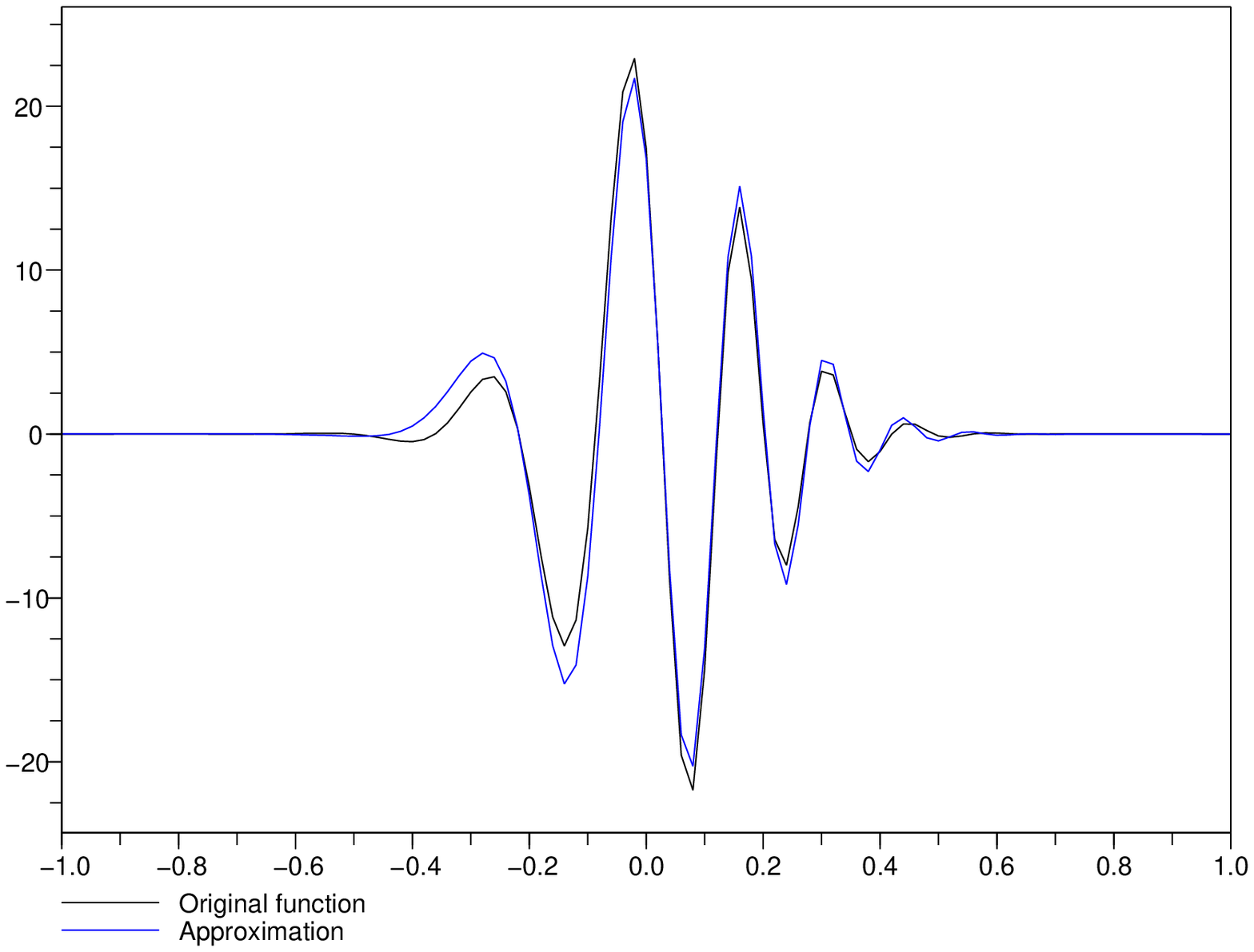,angle=0,width=2.5in,height=3.75in}}
\centerline{\epsfig{file=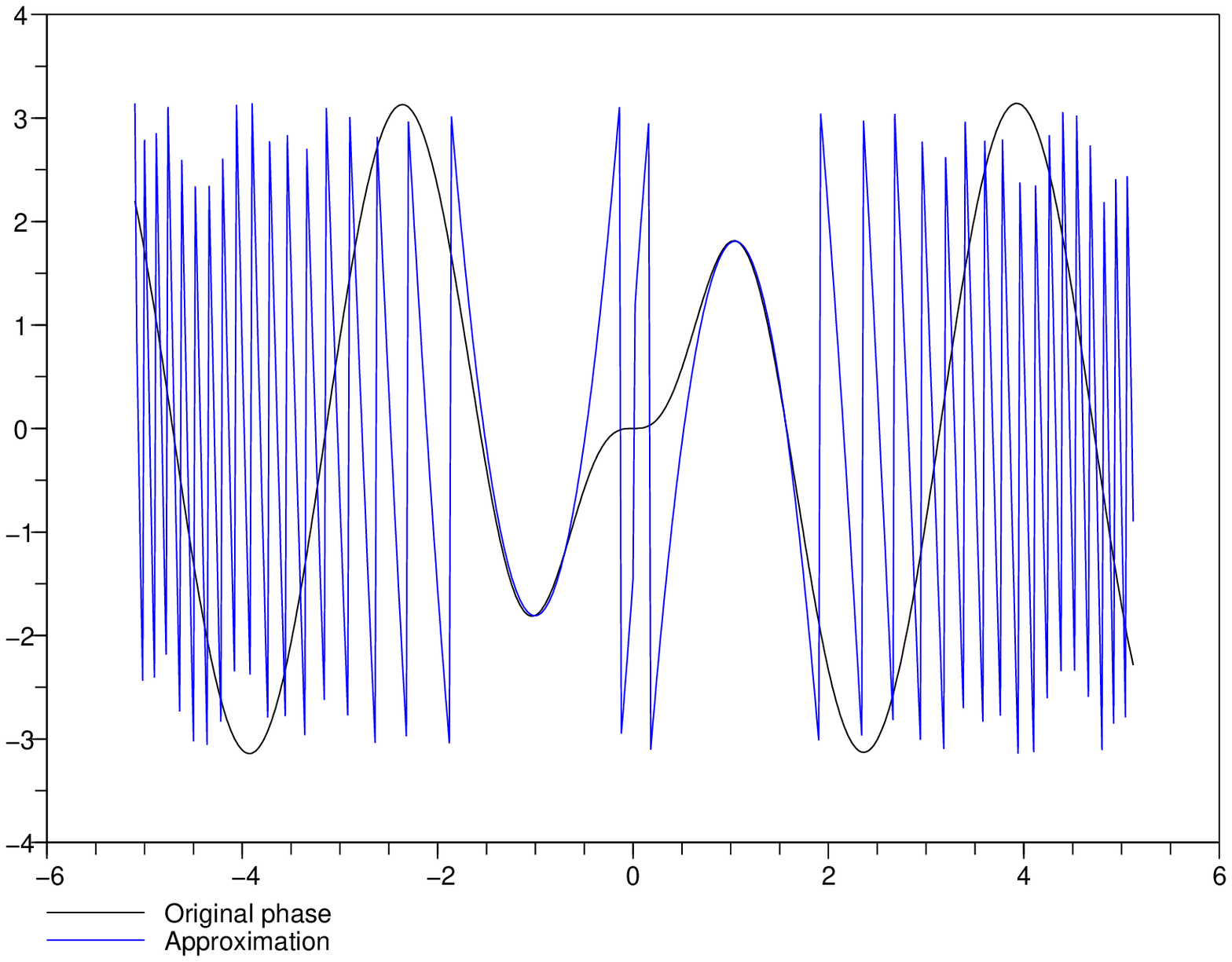,angle=0,width=2.5in,height=3.75in}
\epsfig{file=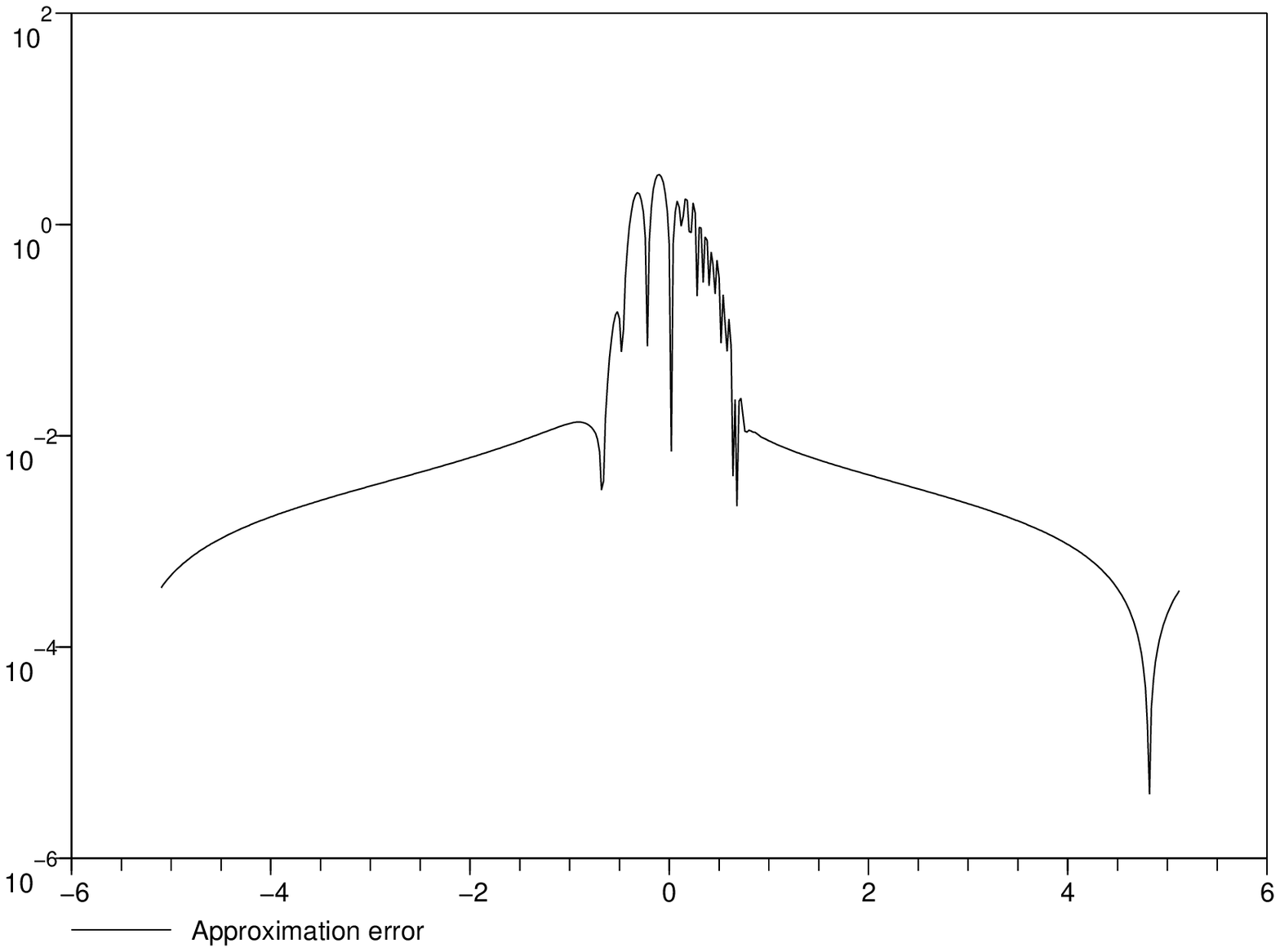,angle=0,width=2.5in,height=3.75in}}
\caption{Pointwise procedure selection on example (\ref{lolo}) for sinusoidal phase.}
\label{sinus}
\end{figure}
\begin{figure}[ht]
\centerline{\epsfig{file=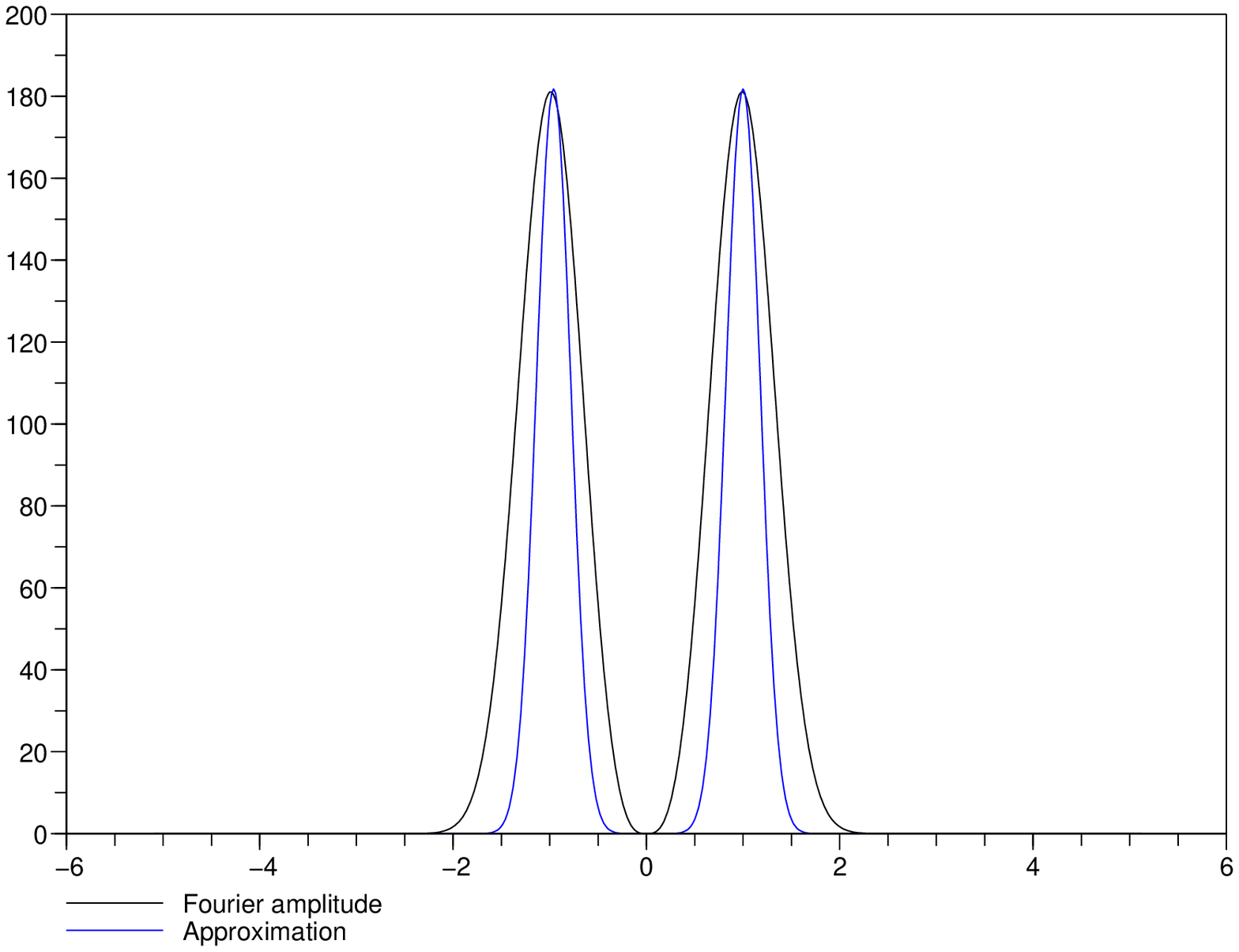,angle=0,width=2.5in,height=3.75in}
\epsfig{file=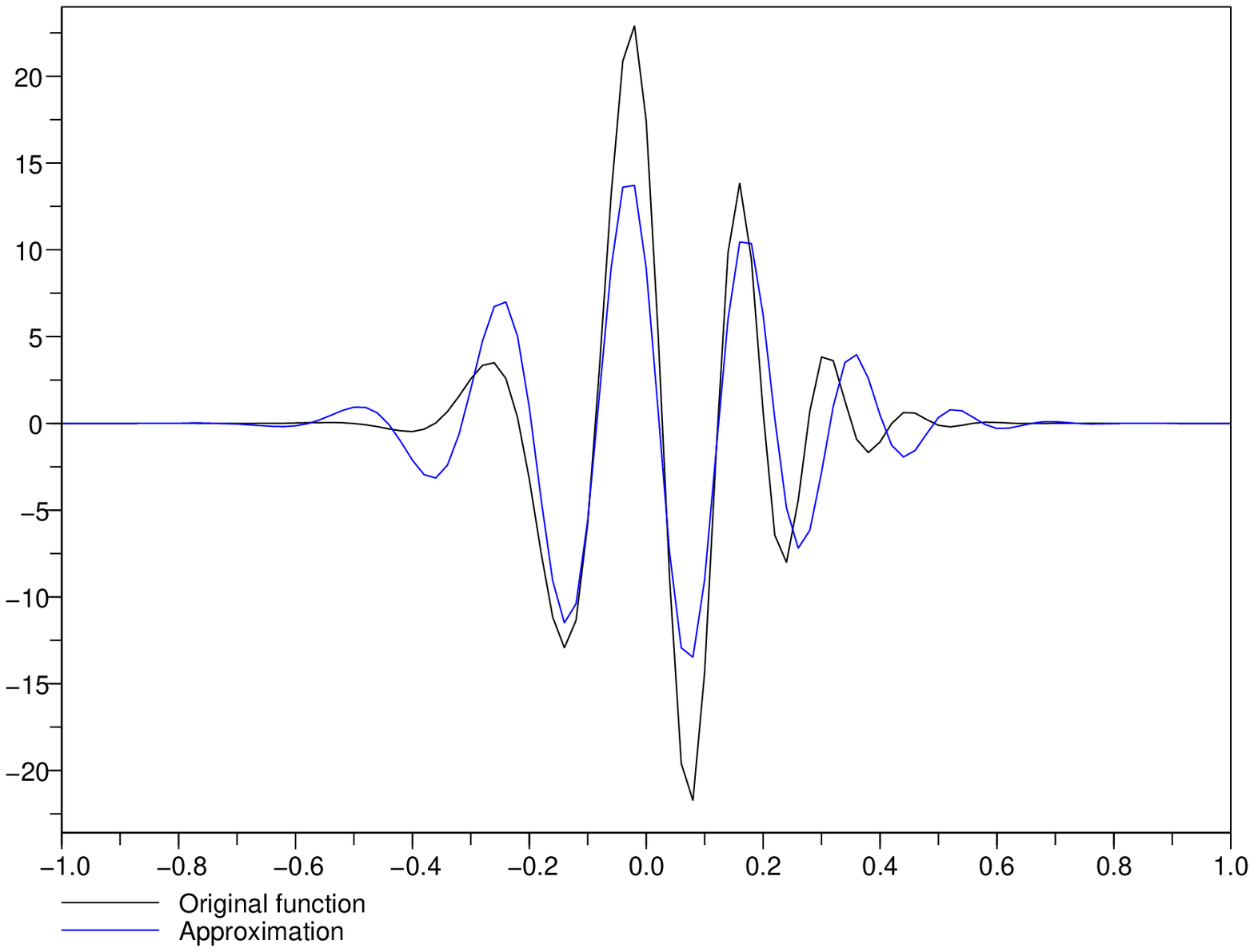,angle=0,width=2.5in,height=3.75in}}
\centerline{\epsfig{file=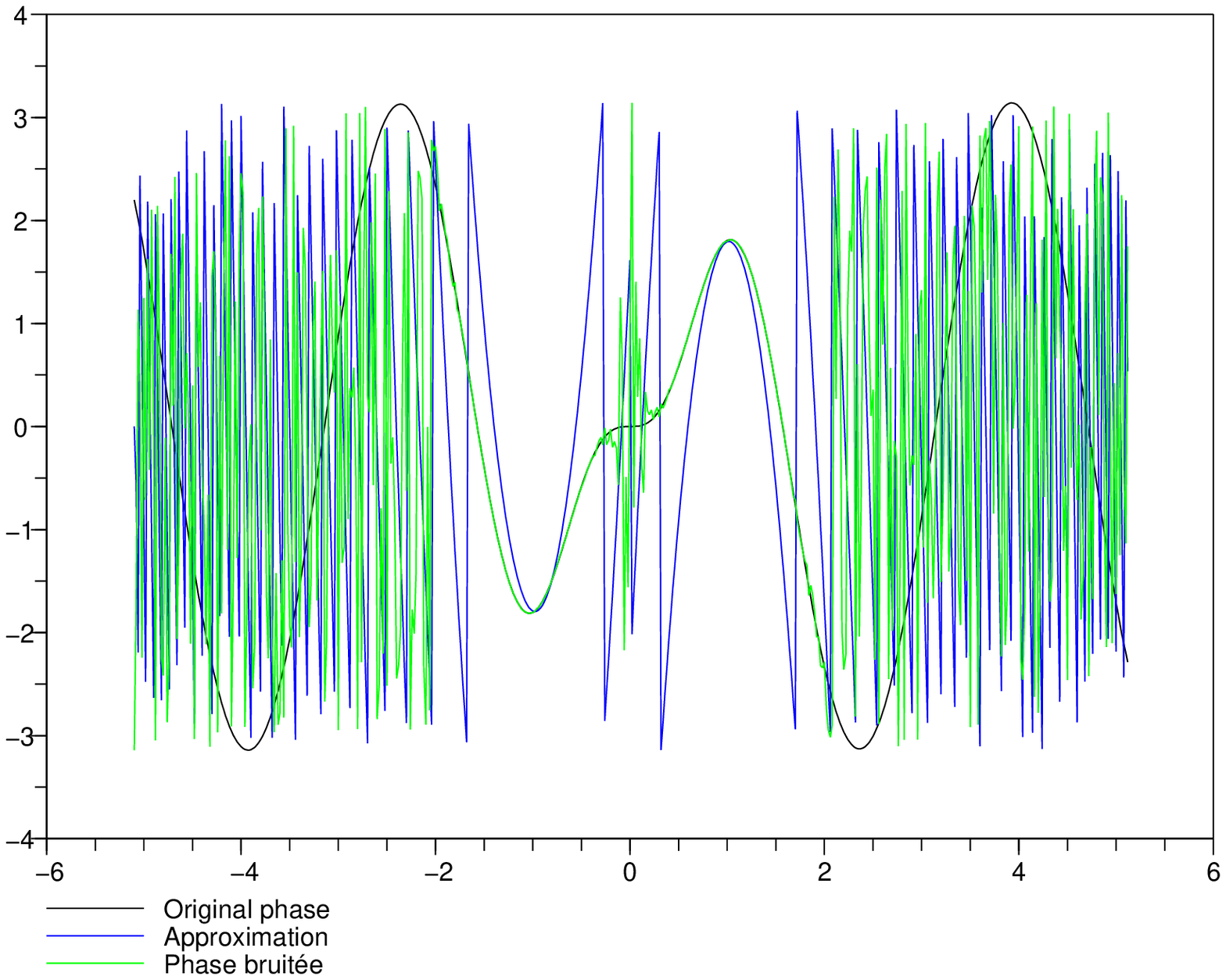,angle=0,width=2.5in,height=3.75in}
\epsfig{file=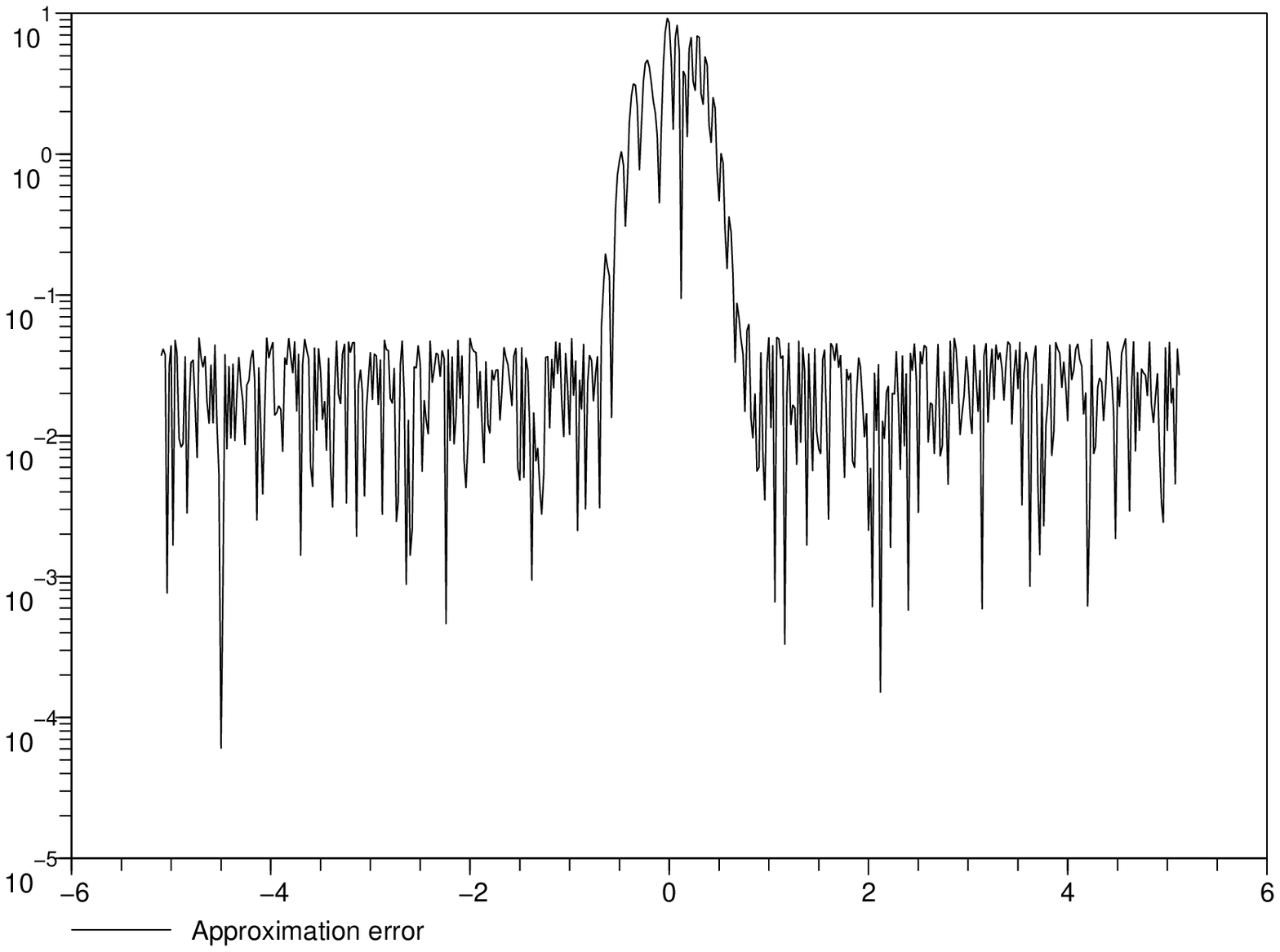,angle=0,width=2.5in,height=3.75in}}
\caption{Pointwise procedure selection on example (\ref{lolo}) slightly corrupted by white noise with sinusoidal phase.}
\label{noise}
\end{figure}
%
%

The specificity of these runs is that now, we test the algorithm ``blind", which means that we just furnish the collection of sample data for both amplitude and phase, and we look after the maximum numerically. Similarly, the derivatives involved in the parameters calculation are computed with finite differences. Corresponding results are displayed in the Figures \ref{cubic}--\ref{noise}. On the left column, one sees the original amplitudes and phases (in black) and their approximations (in blue); on the right one, there are the signal as a function of $t$ (in black), its chirplet approximation (in blue) and finally the absolute error in log-scale. The recovery of the amplitude (\ref{lolo}) by means of two Gaussians looks very satisfying and errors are noticeable only in the last example where the original signal is corrupted by {\it white noise}; the cubic phase (Fig. \ref{cubic}) is approximated in a correct way by second degree polynomials around the maxima of $A$ which are well identified. Of course, in case $A$ would display several local extrema, a more involved process would have to be set up in order to localize them properly inside the vector of samples. Concerning the sinusoidal phase model (Fig. \ref{sinus}), its recovery by means of polynomials of degree two is of course rather poor; however, what really matters is that it is correct in the vicinity of maximum points of $A$, and this is just what happens (see \cite{phase}). The behavior of these approximations in the $t$ variable is good, and the absolute error remains reasonable.

Fig. \ref{noise} deals with a more difficult test-case, namely we set up the amplitude (\ref{lolo}) with the sinusoidal phase, we perform the inverse Fourier transform to get it as a function of $t$, we corrupt the resulting signal with white noise (which makes it not band-limited anymore). This is the data we furnish to the algorithm of \S3.1. Obviously, the white noise has to remain small in order not to perturb the maxima of $A$. What we found is that even this requirement isn't enough as the phase function should not be corrupted too much in order to let the algorithm approximate it locally by a parabola correctly. In this case, the absolute recovery error is bigger compared to the noise-free case (see again Fig. \ref{sinus}).

\subsection{Finding chirp patterns in stock market indices} 

The Hang-Seng Composite Enterprise Index is one of the main index for the Chinese stock market and the CAC 40 is the leading index for the Parisian place. We used the corresponding daily price fluctuations on an interval of 1024 days (ending on august 29th 2008) in order to check on a real-life case whether or not the level of tolerance to noise was acceptable. Clearly, in order to avoid as much as possible spurious local extrema in the Fourier spectrum of the data, we {\it detrended} the prices by a global least-squares interpolation\footnote{This construction ensures the fluctuation has some amount of vanishing moments.}; polynomials of degree 5 and 6 have been used. Moreover, for the CAC 40 only, the detrended fluctuation was displaying quite a sharp peak corresponding to a periodic cycle. We also used a Basis Pursuit algorithm to remove as much as possible the white noise components of these fluctuations, see \cite{donoho}.
\begin{figure}[ht]
\centerline{\epsfig{file=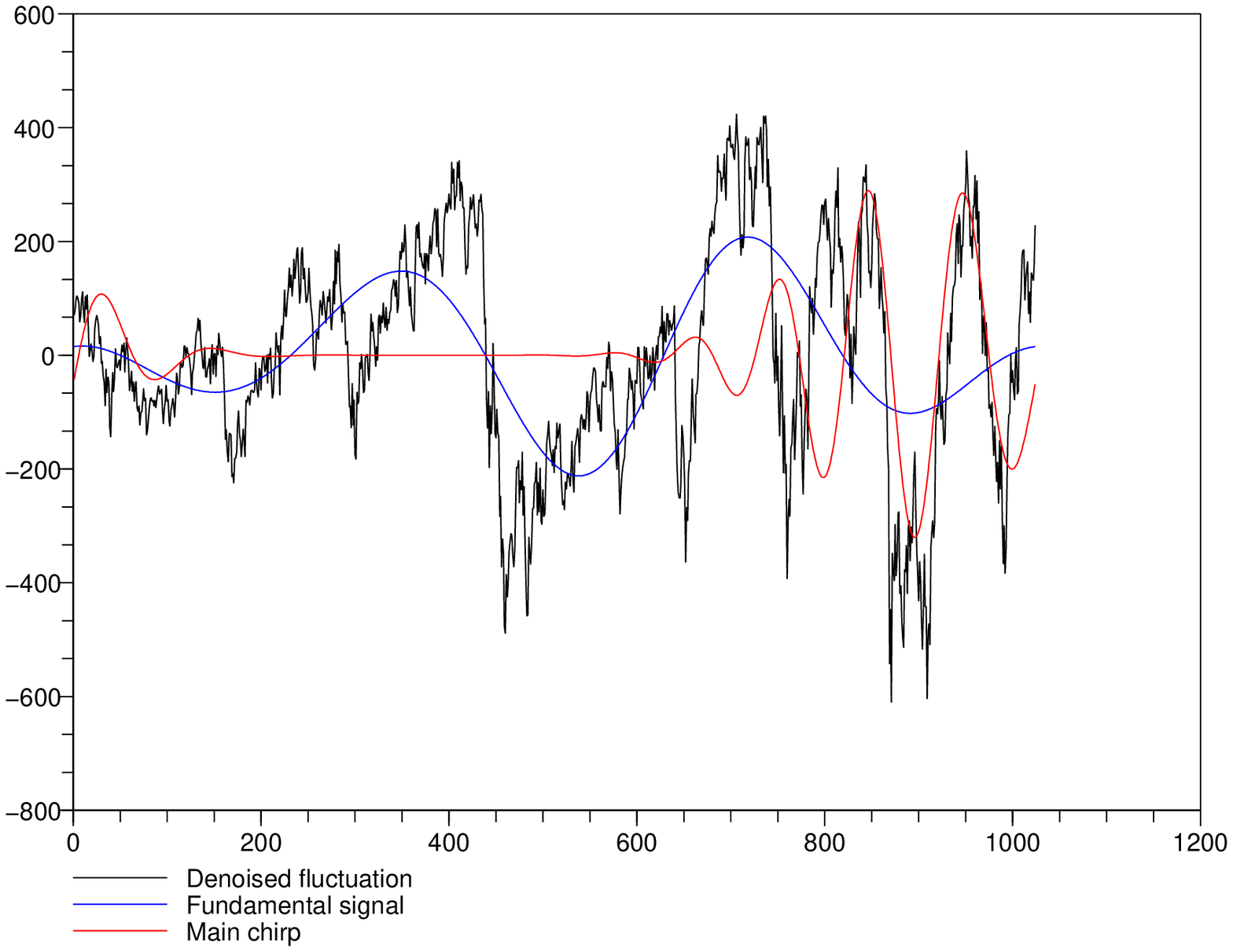,angle=0,width=4.5in,height=3.75in}}
\centerline{\epsfig{file=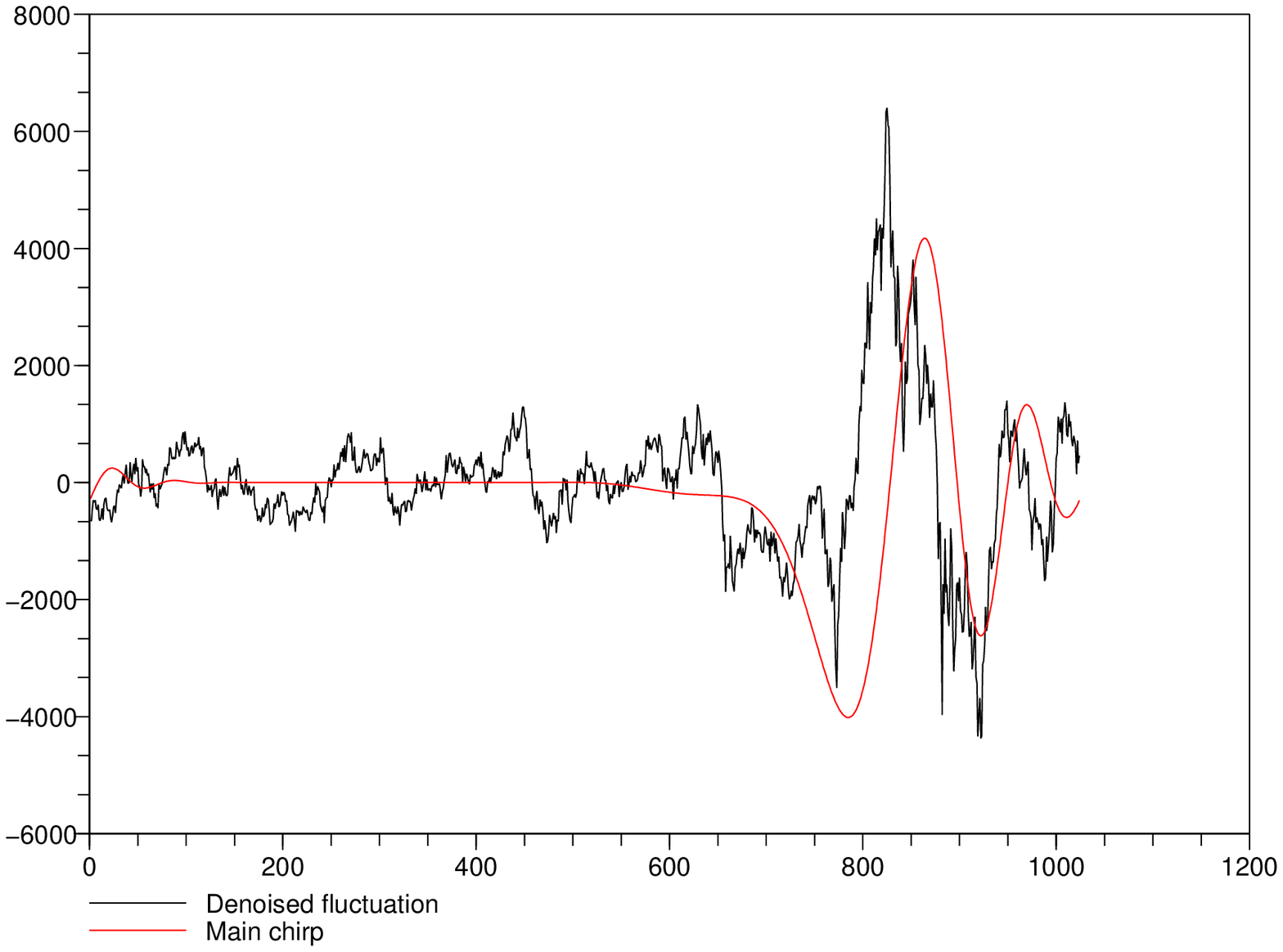,angle=0,width=4.5in,height=3.75in}}
\caption{Chirp patterns in the daily price fluctuations of CAC 40 (top) and HSCEI (bottom).}
\label{hscei}
\end{figure}
On Fig. \ref{hscei}, we display the detrended data (in black), the long-lasting periodic cycle of the CAC 40 (in blue) and the chirps we got out of the algorithm of \S3.1. The case of the Chinese index seem to be quite interesting as the recovery seems to be quite sharp. Results on both CAC 40 and HSCEI suggest that a qualitative change in the behavior of blue chips quotations occurred after the ignition of the so-called {\it subprime crisis/credit crunch} (corresponding to abscissa $t \simeq 700$).

\end{document}